%% file: main.tex
\documentclass [11pt,oneside]{article}

\reversemarginpar \textwidth=14cm \textheight=19cm

\usepackage{amssymb} \usepackage{amsfonts} \usepackage{amsmath}
\usepackage{amsthm} \usepackage{epsfig}
\usepackage{graphics}

\include{defs}

\author{Roberto \titsc{Frigerio}\footnote{Supported by the INTAS project
``CalcoMet-GT'' 03-51-3663}}

\title{Commensurability 
of hyperbolic manifolds\\ with geodesic boundary}

\begin{document}

\maketitle

\begin{abstract}
Suppose $n\geqslant 3$, let $M_1, M_2$ be $n$-dimensional connected 
complete finite-volume 
hyperbolic manifolds with non-empty geodesic boundary, and 
suppose that $\pi_1 (M_1)$ is quasi-isometric to $\pi_1 (M_2)$
(with respect to the word metric).
Also suppose that if $n=3$, then 
$\partial M_1$ and $\partial M_2$ are compact.
We show that $M_1$ is commensurable with $M_2$.
Moreover, we show that there exist
homotopically equivalent hyperbolic $3$-manifolds with non-compact
geodesic boundary which are not commensurable with each other.

We also prove that if $M$ 
is as $M_1$ above and $G$ is a finitely generated group
which is quasi-isometric to $\pi_1 (M)$, then there exists a
hyperbolic manifold with geodesic boundary $M'$ 
with the following properties: $M'$ is commensurable with $M$,
and $G$ is a finite extension of a group which contains
$\pi_1 (M')$ as a finite-index subgroup.
\vspace{4pt}


\vspace{4pt}

\end{abstract}

\noindent
A \emph{quasi-isometry} between metric spaces is a (not necessarily
continuous) map which
looks biLipschitz from afar (see Section~\ref{prel:section} 
for a precise definition). 
Even if they do not give information about
the local structure of metric spaces, quasi-isometries
usually capture the most important properties of 
their large-scale geometry. 

Let $\Gamma$ be any group with a finite set $S$ of generators.
Then a locally finite graph is defined which depends on
$\Gamma$ and $S$ (such a graph is called
a \emph{Cayley graph} of $\Gamma$, see Section~\ref{prel:section}). 
This graph is naturally
endowed with a path metric, and metric graphs arising from
different sets of generators are quasi-isometric to each
other. Thus we can associate to any finitely generated group
a well-defined quasi-isometry class
of locally finite metric graphs. A lot of energies have been devoted
in the last years in order to understand how algebraic and geometric
properties
of a finitely generated group can be related to
metric properties of its Cayley graph (see~\cite{Gro} for
a detailed survey on this topic). This approach
has proved particularly powerful, for example, 
in providing intrinsic 
characterization of fundamental groups of metric spaces with
curvature bounds. 

It is easily seen that (the Cayley graphs of)
the fundamental groups of two \emph{compact} hyperbolic $n$-manifolds
without boundary are always quasi-isometric to each other.
On the other hand, Schwartz proved in~\cite{Sch} that
in dimension greater than $2$, the fundamental groups of
two complete finite-volume non-compact hyperbolic
manifolds without boundary are quasi-isometric to each other
if and only if the manifolds themselves are commensurable with each other
(see below for the definition of commensurability).
In this paper we extend Schwartz's result to the case of hyperbolic
manifolds with non-empty geodesic boundary. Namely, let
$n\geqslant 3$, let $M_1, M_2$ be $n$-dimensional connected 
complete finite-volume 
hyperbolic manifolds with non-empty geodesic boundary, and 
suppose that if $n=3$, then $\partial M_1$ and $\partial M_2$ are compact.
We prove that if
$\pi_1 (M_1)$ is quasi-isometric to $\pi_1 (M_2)$, then
$M_1$ is commensurable with $M_2$. Moreover, we show that there exist
homotopically equivalent hyperbolic $3$-manifolds with non-compact
geodesic boundary which are not commensurable with each other.

We also prove
the following quasi-isometric rigidity result for fundamental groups
of hyperbolic manifolds with geodesic boundary.
Let
$n\geqslant 3$, let $M$ be an $n$-dimensional connected 
complete finite-volume 
hyperbolic manifold with non-empty geodesic boundary, and 
suppose that if $n=3$, then $\partial M$ is compact.
We show that if $G$ is a finitely generated group
which is quasi-isometric to $\pi_1 (M)$, then there exists a
hyperbolic manifold with geodesic boundary $M'$ 
with the following properties: $M'$ is commensurable with $M$,
and $G$ is a finite extension of a group which contains
$\pi_1 (M')$ as a finite-index subgroup.
A similar quasi-isometric rigidity result 
for fundamental groups of complete finite-volume non-compact
hyperbolic manifolds without boundary
was proved by Schwartz in~\cite{Sch}.

\section{Preliminaries and statements}\label{prel:section}
  
In this section we
list some preliminary facts about hyperbolic manifolds and
Cayley graphs of finitely generated groups,
and we state our main results.

\paragraph{Quasi-isometries}
Let $(X_1,d_1)$, $(X_2,d_2)$ be metric spaces and let 
$k\geqslant 1$ be a real number. A map $f:X_1\to X_2$ 
is a \emph{$k$-quasi-isometric embedding} if for every
$x,x'\in X_1$ we have
\[
k^{-1} \cdot d_1 (x,x')-k\leqslant d_2 (f(x),f(x'))
\leqslant k \cdot d_1 (x,x')+k.
\]
The map $f$ is said to be a $k$-\emph{quasi-isometry}
if there exists a $k$-quasi-isometric embedding
$g:X_2\to X_1$ such that
\[
d_1 (x,g(f(x)))\leqslant k,\ d_2 (x',f(g(x')))\leqslant k
\quad \textrm{for every}\ x\in X_1,\, x'\in X_2.
\]
In this case we say that $g$ is a $k$-\emph{pseudo-inverse}
of $f$, and that $(X_1,d_1)$ and $(X_2,d_2)$ are \emph{$k$-quasi-isometric}
to each other. A map is a quasi-isometry if it is a $k$-quasi-isometry
for some $k\in\matR$, and two spaces are quasi-isometric if they
are $k$-quasi-isometric for some $k\in\matR$.

Two quasi-isometries $f,g:X_1\to X_2$ are \emph{equivalent}
if a constant $c>0$ exists such that $d_2 (f (x),g(x))\leqslant c$
for every $x\in X_1$. Equivalence classes
of quasi-isometries of a fixed metric space $X$ form the group
${\rm QIsom} (X)$,
which is called the quasi-isometry group of $X$.

\paragraph{The Cayley graph of a group}
Let $\Gamma$ be a group and $S\subset\Gamma$ be a finite set of
generators for $\Gamma$. We also suppose that $S$ is symmetric,
\emph{i.e.} that $S=S^{-1}$, and that $1\notin S$. 
For any $\gamma\in\Gamma$ let
$|\gamma|_S$ be the minimal length of words 
representing $\gamma$ and having letters in $S$. 
If $\gamma_1,\gamma_2$ are elements of $\Gamma$  
we set $d_S (\gamma_1,\gamma_2)=|\gamma_2^{-1}\gamma_1|_S$. 
It is easily seen that $d_S$ is a distance which turns $\Gamma$
into a discrete metric space. Moreover, left translations
by elements in $\Gamma$ act as isometries with respect to $d_S$.
We are now ready to define the \emph{Cayley graph}
$\calC (\Gamma,S)$ of $\Gamma$ (with respect to $S$) as follows:
the vertices of 
$\calC (\Gamma,S)$
are the elements of $\Gamma$, two vertices $\gamma_1,\gamma_2\in\Gamma$
are joined by an edge if and only if $d_S (\gamma_1,\gamma_2)=1$
(\emph{i.e.}~if and only if $\gamma_1=\gamma_2 s$ for some $s\in S$),
and $\calC (\Gamma,S)$ does not contain multiple edges
(\emph{i.e.}~two vertices are connected by at most one edge).
It is easily seen that $d_S$ can be extended to a distance (still denoted
by $d_S$) on $\calC (\Gamma,S)$ which turns $\calC (\Gamma,S)$
into a path metric space. Moreover, the action of $\Gamma$ 
onto itself by left translations can be extended to an isometric
action of $\Gamma$ on $\calC (\Gamma,S)$.

It is easily seen that if $S,S'$ are finite symmetric sets of generators
for $\Gamma$, then $\calC (\Gamma,S)$ is quasi-isometric to
$\calC (\Gamma,S')$. We say that 
a finitely generated group $\Gamma$ is quasi-isometric to
a metric space $(X,d)$
if some (and then any) Cayley graph of $\Gamma$ is quasi-isometric 
to $(X,d)$. In particular, 
two finitely generated
groups $\Gamma_1,\Gamma_2$ are quasi-isometric if some (and then any)
Cayley graph of $\Gamma_1$ is quasi-isometric to some (and then any)
Cayley graph of $\Gamma_2$. 
Every finite group is quasi-isometric
to the trivial group. More in general, it is easily seen that 
every finite-index subgroup of any given group is quasi-isometric
to the group itself. In a similar way, every finite extension of any
given group is quasi-isometric to the group itself.

If $\Gamma'$ is a subgroup of $\Gamma$, we denote by $\Gamma/\Gamma'$
the set of left lateral classes of $\Gamma'$ in $\Gamma$, and by
$[\Gamma : \Gamma']=\# (\Gamma/\Gamma')$ the index of $\Gamma'$
in $\Gamma$.

\paragraph{An important relationship with the universal covering}
The (easy) result we are going to describe~\cite{Mil,Ghy}
establishes perhaps the most 
important
relationship between
the geometry of a compact metric space and the metric
properties of (a Cayley graph of) its fundamental group. 

Let $(X,d)$ be a compact Riemannian manifold, possibly with boundary.
Let
$(\widetilde{X},\widetilde{d})$ be the metric universal covering
of $(X,d)$, let $\widetilde{x}\in\widetilde{X}$
be a given basepoint and 
fix an identification of $\pi_1 (X)$ with the group of the covering
automorphisms of $\widetilde{X}$. Then the map 
\[
\pi_1 (X)\ni \gamma \mapsto \gamma (\widetilde{x})\in\widetilde{X}
\]
is a $k$-quasi-isometry between $\pi_1 (X)$ and 
$(\widetilde{X},\widetilde{d})$, where $k$ only depends
on the diameter of $X$ and
the set of generators used to define the Cayley graph
of $\pi_1 (X)$.

Since the universal covering of any compact hyperbolic
$n$-manifold without boundary is isometric to $\matH^n$,
this implies that the fundamental groups
of any two compact hyperbolic 
$n$-manifolds without boundary are quasi-isometric to each
other.
On the contrary,
Theorem~\ref{main1:teo} below implies that there exist pairs of
compact
hyperbolic manifolds
with non-empty geodesic boundary having non-quasi-isometric 
fundamental groups.

\paragraph{Commensurability of hyperbolic manifolds}

From now on we will always suppose $n\geqslant 3$.
Moreover, all manifolds will be connected.
By \emph{hyperbolic $n$-manifold} we will mean 
a complete finite-volume 
hyperbolic $n$-manifold with (possibly empty)
geodesic boundary. Let $N_1,N_2$ be hyperbolic $n$-manifolds.
We say that $N_1$ is \emph{commensurable}
with $N_2$ if a hyperbolic $n$-manifold $N_3$ exists which is 
the total space of a finite Riemannian covering both of $N_1$
and of $N_2$. Notice that in this case the fundamental group 
of $N_3$ is a finite-index subgroup both of $\pi_1 (N_1)$
and of $\pi_1 (N_2)$, so 
$\pi_1 (N_1)$ is quasi-isometric to
$\pi_1 (N_2)$. The main result of this paper shows that the converse
of this statement is also true, under some additional conditions:

\begin{teo}\label{main1:teo}
Let $N_1,N_2$ be hyperbolic $n$-manifolds with non-empty geodesic
boundary and assume that if $n=3$ then $\partial N_1$ and
$\partial N_2$ are compact.
Suppose that $\pi_1 (N_1)$ is quasi-isometric
to $\pi_1 (N_2)$. 
Then $N_1$ is commensurable with $N_2$.
\end{teo}

\paragraph{The fundamental group as a group of isometries}
We denote by $(\matH^n,d^\matH)$ (or simply by $\matH^n$)
the hyperbolic $n$-space.
Let $N$ be a hyperbolic $n$-manifold with non-empty
boundary
and let $\pi:\widetilde{N}\to N$ be
the universal covering of $N$. By developing $\widetilde{N}$
in $\matH^n$ we can identify $\widetilde{N}$ with
a convex polyhedron of $\matH^n$ bounded by
a countable number of disjoint geodesic hyperplanes.
The group of the automorphisms of the covering
$\pi:\widetilde{N}\to N$ can be identified in a natural way
with 
a discrete torsion-free subgroup $\Gamma$ of 
$\mathrm{Isom} (\widetilde{N})<\mathrm{Isom}(\matH^n)$ such that
$N\cong\widetilde{N}/\Gamma$.
Also recall that there exists an isomorphism
$\pi_1(N)\cong\Gamma$, which is canonical up to
conjugacy. With a slight abuse, from now on we will often refer to 
$\Gamma$ as to the fundamental
group of $N$. Notice that $N$ uniquely determines $\Gamma$ up
to conjugation by elements in ${\rm Isom}(\matH^n)$.

\paragraph{The commensurator}
Let $\Gamma_1,\Gamma_2$ be discrete subgroups of
${\rm Isom} (\matH^n)$. We say that an element $g\in{\rm Isom}
(\matH^n)$ \emph{commensurates} $\Gamma_1$ with $\Gamma_2$
if $g\Gamma_1 g^{-1}\cap \Gamma_2$ has finite index both
in $g\Gamma_1 g^{-1}$ and in $\Gamma_2$. 
In this case we also say that $\Gamma_1$ and $\Gamma_2$ are commensurable
with each other. We will see in Lemma~\ref{commensurators:lemma}
that two hyperbolic manifolds with non-empty geodesic boundary
$N_1,N_2$ are commensurable if and only if their fundamental
groups are.

If $\Gamma$ is a discrete subgroup of ${\rm Isom} (\matH^n)$,
the group of elements of ${\rm Isom} (\matH^n)$ which commensurate
$\Gamma$ with itself is called the \emph{commensurator}
of $\Gamma$, and it is denoted by ${\rm Comm} (\Gamma)$.
Of course $\Gamma$ itself is a (not necessarily normal) subgroup
of ${\rm Comm} (\Gamma)$.

\paragraph{Commensurator and quasi-isometry group}
Let $G$ be a finitely generated group and $g$ be an element of $G$.
We denote by $\ell_G (g)\in{\rm QIsom} (G)$ the
equivalence class of the  
quasi-isometry induced by the left translation by $g$. 
The map $\ell_G :G \to
{\rm QIsom} (G)$ is a homomorphism,
whose image is not necessarily normal in
${\rm QIsom} (G)$. 
In what follows we will often denote
just by $g$ the element $\ell_G (g)\in{\rm QIsom} (G)$,
and we will consider $G$ as a subgroup of ${\rm QIsom} (G)$
via the homomorphism $\ell_G$. (In the cases we will be concerned with, 
the homomorphism $\ell_G$ is injective. Note however that
this is not true in general: for example, if $G$ is Abelian
then $\ell_G$ is the trivial homomorphism.) 
The following result will be proved in Section~\ref{corollari:section}.

\begin{teo}\label{comm1:teo}
Let $N$ be a hyperbolic $n$-manifold with non-empty boundary,
let $\Gamma<{\rm Isom} (\matH^n)$ be the fundamental group of $N$
and if $n=3$ also suppose that $\partial N$ is compact.
Then the identity of $\Gamma$ extends to a
canonical isomorphism between
the commensurator of $\Gamma$
and the quasi-isometry group of $\Gamma$.
In particular, we have 
$[{\rm Comm} (\Gamma):\Gamma]=
[{\rm QIsom} (\Gamma):\Gamma]$.  
\end{teo}

If $G<{\rm Isom} (\matH^n)$ 
is the fundamental group of a non-compact complete finite-volume
hyperbolic manifold without boundary, results of Margulis and 
Borel~\cite{Bor,Zim}
imply that ${\rm Comm} (G)/G$ is finite if and only if $G$ is non-arithmetic.
Moreover, a theorem of Schwartz~\cite{Sch} ensures that 
$[{\rm Comm} (G):G]=[{\rm QIsom} (G):G]$, 
so arithmeticity of $G$
turns out to depend only on the quasi-isometry type of $G$.
The following result will be proved at the end of 
Section~\ref{ciccia:section}, and shows that things look very different
in the case of non-empty geodesic boundary.

\begin{prop}\label{comm2:teo}
Let $N$ be a hyperbolic $n$-manifold with non-empty boundary and 
let $\Gamma<{\rm Isom} (\matH^n)$ be the fundamental group of $N$. Then
${\rm Comm} (\Gamma)/\Gamma$ is finite. 
\end{prop}

Thus, if 
$n\geqslant 4$ or $\partial N$ is compact,
then ${\rm QIsom} (\pi_1 (N))/\pi_1 (N)$
is a finite group.

\paragraph{Quasi-isometric rigidity}
In Section~\ref{corollari:section} we will prove that
quasi-isometries essentially preserve the property of being the
fundamental group of a hyperbolic manifold with boundary:

\begin{teo}\label{main3:teo}
Let $N$ be a hyperbolic $n$-manifold with non-empty geodesic boundary,
and if $n=3$ also suppose that $\partial N$ is compact. 
Let $\Gamma<{\rm Isom}(\matH^n)$ be the fundamental group
of $N$ and $G$ be a finitely generated
abstract group which is quasi-isometric to $\Gamma$.
Then $G$ is a finite extension of a discrete subgroup 
$\Gamma'<{\rm Isom} (\matH^n)$ which is commensurable with $\Gamma$.
\end{teo}

By Selberg Lemma~\cite{Sel}, any finitely generated discrete 
subgroup of ${\rm Isom} (\matH^n)$ contains
a finite-index torsion-free subgroup, so Theorem~\ref{main3:teo}
implies the following:

\begin{cor}\label{main4:teo}
Let $N$ and $G$ be as in the statement of Theorem~\ref{main3:teo}.
Then there exists a hyperbolic manifold $N'$ with the following
properties: 
$N'$ is commensurable with $N$, and
$G$ is a finite extension of a group which contains $\pi_1 (N')$
as a finite-index subgroup.
\end{cor}

\paragraph{Counterexamples}
In Section~\ref{contro:section} 
we will show that the hypotheses of Theorems~\ref{main1:teo},
\ref{comm1:teo} cannot be
weakened. Namely,
we will prove the following:

\begin{teo}\label{contro1:teo}
There exist non-commensurable hyperbolic $3$-manifolds with 
non-compact geodesic boundary 
sharing the same fundamental group.
\end{teo}

\begin{teo}\label{contro2:teo}
A hyperbolic $3$-manifold $M$ with non-compact geodesic boundary
exists such that $\pi_1 (M)$
has infinite index in ${\rm QIsom} (\pi_1 (M))$.
\end{teo}

\section{Hyperbolic manifolds with geodesic boundary}\label{ciccia:section}

This section is devoted to a brief description of the most important
topological and geometric properties of hyperbolic
manifolds with geodesic boundary.

\paragraph{Natural compactification}  
Let $N$ be a hyperbolic $n$-manifold with non-empty geodesic boundary.
Then $\partial N$, endowed
with the Riemannian metric it inherits from $N$, is a
hyperbolic $(n-1)$-manifold without boundary
(completeness of $\partial N$ is obvious, and the 
volume of $\partial N$ is proved to be finite in~\cite{Kojima1}).
It is well-known~\cite{Kojima1, Kojima2} that $N$ consists of a compact portion
together with some cusps based on Euclidean $(n-1)$-manifolds. 
More precisely, the $\varepsilon$-thin part of $N$
(see~\cite{thu}) consists of cusps of the form $T\times [0,\infty)$,
where $T$ is a compact Euclidean $(n-1)$-manifold 
with (possibly empty) geodesic boundary such that 
$\left(T\times[0,\infty)\right)\cap\partial N=
\partial T\times[0,\infty)$.
A cusp based on a closed Euclidean $(n-1)$-manifold
is therefore disjoint from $\partial N$ and
is called \emph{internal}, while a cusp
based on a Euclidean $(n-1)$-manifold with non-empty
boundary intersects $\partial N$ in one or two
internal cusps of $\partial N$, 
and is called a \emph{boundary cusp}.  
This
description of the ends of $N$ easily implies that $N$ admits
a natural compactification $\overline{N}$ obtained by adding
a closed Euclidean $(n-1)$-manifold for each internal cusp
and a compact Euclidean $(n-1)$-manifold with non-empty 
geodesic boundary
for each boundary cusp. 

For later purposes, we observe that when $n=3$, 
$\overline{N}$ is obtained by adding to $N$
some tori, Klein bottles, closed annuli and closed Moebius strips.

\paragraph{Universal covering}
Let $\pi:\widetilde{N}\to N$ be
the universal covering of $N$. We recall that $\widetilde{N}$
can be identified with
a convex polyhedron of $\matH^n$ bounded by
a countable number of disjoint geodesic hyperplanes $S^i,\ i\in\matN$.
For any $i\in\matN$ let $S^i_+$ denote the closed 
half-space of $\matH^n$ 
bounded by $S^i$ and containing $\widetilde{N}$, let
$S^i_-$ be the closed half-space of $\matH^n$ opposite to
$S^i_+$ and let $\Delta^{\! i}$ be the internal part of the closure
at infinity of $S^i_-$. Of course we have $\widetilde{N}=
\bigcap_{i\in\matN} S^i_+$, so denoting by 
$\widetilde{N}_{\infty}$ the closure at infinity of $\widetilde{N}$
we obtain $\widetilde{N}_{\infty}=\partial\matH^n\setminus
\bigcup_{i\in\matN} \Delta^{\! i}$. We also denote by $\overline{S}^{i}$
the closure at infinity of ${S}^{i}$. 

It is easily seen that any internal cusp of $N$ lifts in $\widetilde{N}$ 
to the union of a countable number of disjoint
horoballs, each of which is entirely contained in $\widetilde{N}$.
Things become a bit more complicated when considering boundary cusps.
First observe that if $i\neq j$, then
$\overline{S}^i\cap\overline{S}^j$ is either empty
or consists of one point in $\partial \matH^n$. 
It is shown in~\cite{Kojima1} that if $q\in
\overline{S}^i\cap\overline{S}^j$, then
the intersection of $\widetilde{N}$ with a sufficiently small
horoball centered at $q$ projects onto a boundary cusp of $N$.
Conversely,
any component of the preimage of 
a boundary cusp of $N$ is the intersection
of $\widetilde{N}$ with a horoball centered at a
point which belongs to the boundary at infinity
of two different components of $\partial\widetilde{N}$.

\paragraph{Limit set and discreteness of ${\rm Isom}(\widetilde{N})$} 
Let $\Gamma<{\rm Isom}(\widetilde{N})<\mathrm{Isom}(\matH^n)$ 
be the fundamental group of $N$,
let $\Lambda(\Gamma)$ denote the limit set of $\Gamma$
and set $\Omega(\Gamma)=\partial \matH^n\setminus\Lambda(\Gamma)$.
Kojima has shown in~\cite{Kojima1} that $\Lambda(\Gamma)=
\widetilde{N}_{\infty}$, so 
the round balls $\Delta^{\! i},\ i\in\matN$ previously
defined actually are the connected components of $\Omega(\Gamma)$.
Since $\widetilde{N}_{\infty}=\Lambda(\Gamma)$,
we have that $\widetilde{N}$
is the intersection of $\matH^n$ 
with the convex hull of $\Lambda(\Gamma)$, so
$N$ is 
the convex core (see~\cite{thu}) of the hyperbolic manifold
$\matH^n/\Gamma$. Thus $\Gamma$ is geometrically finite
and it uniquely determines
$N$.
In Section~\ref{cicciavera:section} we will need the following result,
which is proved in~\cite{Fri} (see also~\cite{KMS}).

\begin{lemma}\label{connect:lemma}
Let $j:S^{n-2}\to\Lambda(\Gamma)$ be a topological
embedding. Then $\Lambda(\Gamma)\setminus j(S^{n-2})$
is path connected if and only if $j(S^{n-2})=\partial \Delta^{\! l}$ for
some $l\in\matN$.  
\end{lemma}

We also have the following: 

\begin{lemma}\label{discrete:lemma}
The group ${\rm Isom} (\widetilde{N})$ of isometries of 
$\widetilde{N}$ is a discrete subgroup of ${\rm Isom} (\matH^n)$.
\end{lemma}
\noindent\emph{Proof:}
Since $\widetilde{N}$ is connected and has non-empty interior,
every element in ${\rm Isom} (\widetilde{N})$ uniquely
determines an element in ${\rm Isom} (\matH^n)$.
Let now $\{g_n\}_{n\in\matN}\subset{\rm Isom} (\widetilde{N})$ 
be a sequence
converging to the identity of $\matH^n$. If $S^i,S^j$ are
components of $\partial\widetilde{N}$ with
$\overline{S}^{i}\cap\overline{S}^{j}=\emptyset$,
we denote by $e_{ij}\subset \widetilde{N}$ the 
unique shortest path joining $S^i$ with $S^j$ and we set
\[
E=\left(\bigcup_{\overline{S}^i\cap\overline{S}^j=\emptyset} e_{ij}\right).
\]
Now $E$ is a $\Gamma$-invariant subset of $\matH^n$ and
the limit set $\Lambda (\Gamma)$ is not contained in the boundary at infinity
of any hyperbolic hyperplane. Thus
there exist 
components $S^1,\ldots,S^{n+1}$ of $\partial \widetilde{N}$
such that 
$\overline{S}^i \cap\overline{S}^j=\emptyset$ for all $i,j=1,\ldots,n+1$,
$i\neq j$, and
the convex hull of the $e_{ij}$'s, $i,j=1,\ldots,n+1$,
$i\neq j$, has non-empty interior in $\matH^n$.
Since $g_k$ tends to the identity as $k$ tends to $\infty$, there exists
$M\gg 0$ such that 
$g_k (S^i)=S^i$ for all $i=1,\ldots,n+1$, $k\geqslant M$.
Thus if $k\geqslant M$ then
$g_k$ restricts to the identity on each $e_{ij}$,
$i,j=1,\ldots,n+1$,
$i\neq j$, whence on the convex hull of such $e_{ij}$'s, which has
non-empty interior. This forces $g_k={\rm Id}_{\matH^n}$ for every
$k\geqslant M$, and proves that ${\rm Isom}(\widetilde{N})$
is discrete. 
\finedimo

\paragraph{Three lemmas}
We are now ready to prove the following:

\begin{lemma}\label{commensurators:lemma}
Suppose $N_1,N_2$ are hyperbolic $n$-manifolds with non-empty 
geodesic boundary, let $\widetilde{N}_i\subset\matH^n$ 
be the universal covering of $N_i$ and 
$\Gamma_i<{\rm Isom} (\widetilde{N}_i)<{\rm Isom} (\matH^n)$
be the fundamental group of $N_i$. Then $N_1$ is commensurable
with $N_2$ if and only if $\Gamma_1$ is commensurable with $\Gamma_2$. 
\end{lemma}
\noindent\emph{Proof:}
Suppose $N_3$ is a finite Riemannian covering both of $N_1$ and of $N_2$
and let $\Gamma_3<{\rm Isom}(\widetilde{N}_3)<{\rm Isom}(\matH^n)$
be the fundamental group of $N_3$. For $i=1,2$, 
the covering projection $p_i:N_3\to N_i$ induces an isometry
$\widetilde{p}_i:\widetilde{N}_3\to\widetilde{N}_i$. Now conjugation by the
isometry $\widetilde{p}_2\circ\widetilde{p}_1^{-1}$ takes 
$(p_1)_\ast (\Gamma_3)$ onto $(p_2)_\ast (\Gamma_3)$. Since 
$(p_i)_\ast (\Gamma_3)$ has finite index in $\Gamma_i$, this implies that
$\Gamma_1$ is commensurable with $\Gamma_2$.

On the other hand, let $g\in{\rm Isom}(\matH^n)$ be an element such
that $\Gamma_3=g\Gamma_1 g^{-1}\cap \Gamma_2$ has finite
index both in $\Gamma_2$ and in $g\Gamma_1 g^{-1}$.
Then $\Lambda (\Gamma_3)=\Lambda (\Gamma_2)=g(\Lambda (\Gamma_1))$, whence
$\widetilde{N}_2=g(\widetilde{N}_1)$ and 
$N_3=\widetilde{N}_2/\Gamma_3$ is hyperbolic with non-empty boundary
(and with universal covering $\widetilde{N}_3=\widetilde{N}_2$).
Of course the natural projection $N_3\cong \widetilde{N}_2/\Gamma_3\to
\widetilde{N}_2/\Gamma_2=N_2$ is a finite Riemannian covering, and
$g^{-1}:
\widetilde{N}_3\to\widetilde{N}_1$ induces a finite Riemannian covering
$\widetilde{N}_3/\Gamma_3\cong N_3\to N_1\cong \widetilde{N}_1/\Gamma_1$.
Thus
$N_1$ is commensurable with $N_2$.
\finedimo

The following result plays a crucial r\^ole in the proof of 
Theorem~\ref{main1:teo}.

\begin{lemma}\label{important:lemma}
Suppose $N_1,N_2$ are hyperbolic $n$-manifolds with non-empty 
geodesic boundary, let $\widetilde{N}_i\subset\matH^n$ 
be the universal covering of $N_i$ and 
$\Gamma_i<{\rm Isom} (\widetilde{N}_i)<{\rm Isom} (\matH^n)$
be the fundamental group of $N_i$. Suppose that $g\in{\rm Isom}
(\matH^n)$ is such that $g(\widetilde{N}_1)=\widetilde{N}_2$.
Then $g$ commensurates $\Gamma_1$ with $\Gamma_2$. In particular,
$N_1$ is commensurable with $N_2$.
\end{lemma}
\noindent\emph{Proof:}
Let $\Gamma$ be the group of isometries of $\widetilde{N}_1$.
Of course we have $\Gamma_1<\Gamma$ and $\Gamma_2<{\rm Isom} 
(\widetilde{N}_2)=g\Gamma g^{-1}$.
By Lemma~\ref{discrete:lemma}, $\Gamma$ is discrete,
so $\widetilde{N}_1/\Gamma$ is isometric to an orbifold
$N_{\rm orb}$ of positive
volume. Of course we also have $N_{\rm orb}\cong
\widetilde{N}_2/(g\Gamma g^{-1})$, so
\[
\begin{array}{c}
{[ g\Gamma g^{-1}: g\Gamma_1 g^{-1} ]}=
{[ \Gamma:\Gamma_1 ]}={\rm vol}\, N_1/{\rm vol}\, N_{\rm orb}<\infty,\\
{[ g\Gamma g^{-1}: \Gamma_2 ]}={\rm vol}\, N_2/{\rm vol}\, N_{\rm orb}<\infty.
\end{array}
\]
Thus $g \Gamma_1 g^{-1}\cap
\Gamma_2$ has finite index in $g \Gamma g^{-1}$, whence \emph{a fortiori}
$g \Gamma_1 g^{-1}\cap
\Gamma_2$ has finite index both in $g \Gamma_1 g^{-1}$ and in $\Gamma_2$.
\finedimo

The following lemma readily implies Proposition~\ref{comm2:teo}.

\begin{lemma}\label{fare:lemma}
Let $N$ be a hyperbolic $n$-manifold with non-empty 
geodesic boundary, let $\widetilde{N}\subset\matH^n$ 
be the universal covering of $N$ and 
$\Gamma<{\rm Isom} (\widetilde{N})<{\rm Isom} (\matH^n)$
be the fundamental group of $N$. Then ${\rm Comm} (\Gamma)={\rm Isom}
(\widetilde{N})$ and
${\rm Comm}(\Gamma)/\Gamma$ is finite. 
\end{lemma}
\noindent\emph{Proof:}
By Lemma~\ref{discrete:lemma}, ${\rm Isom}(\widetilde{N})$ is discrete,
so $\widetilde{N}/{\rm Isom}(\widetilde{N})$ is isometric to an orbifold
$N_{\rm orb}$ of positive
volume. So
\begin{equation}\label{ciao:eq}
[{\rm Isom}(\widetilde{N}):\Gamma]={\rm vol}\, N/{\rm vol}\, N_{\rm orb}
<\infty. 
\end{equation}
Thus for every $g\in{\rm Isom} (\widetilde{N})$
we have 
\[
[\Gamma:\Gamma\cap g\Gamma g^{-1}]<\infty,\quad 
[g\Gamma g^{-1}:\Gamma\cap g\Gamma g^{-1}]<\infty.
\] 
This implies that
${\rm Isom}(\widetilde{N})$ is contained in ${\rm Comm}(\Gamma)$.

On the other hand, if $g\in{\rm Isom}(\matH^n)$ 
belongs to ${\rm Comm} (\Gamma)$ then 
$\Lambda (\Gamma)=\Lambda (g\Gamma g^{-1})=g (\Lambda (\Gamma))$. Since 
$\widetilde{N}$ is the hyperbolic convex hull of
$\Lambda(\Gamma)$, it follows that $g (\widetilde{N})=\widetilde{N}$,
whence $g\in{\rm Isom}(\widetilde{N})$. Thus ${\rm Comm} (\Gamma)=
{\rm Isom} (\widetilde{N})$, and the conclusion follows 
from inequality~(\ref{ciao:eq}).
\finedimo

\section{Quasi-isometries of hyperbolic polyhedra}\label{cicciavera:section}
Let $N$ be a hyperbolic $n$-manifold with non-empty 
geodesic boundary and denote by $\widetilde{N}\subset\matH^n$ 
the universal covering of $N$.

\paragraph{Neutered hyperbolic polyhedra}
Let $N^\ast\subset N$ be a compact core of $N$ whose
preimage $\widetilde{N^\ast}$ in $\widetilde{N}$ 
is given by the complement in $\widetilde{N}$ of a 
countable family of disjoint horoballs. 
We will refer to $\widetilde{N^\ast}$ as to
the \emph{neutered} universal covering of $N$ (this terminology
is taken from~\cite{Sch}). 
The components of $\widetilde{N}\setminus
\widetilde{N^\ast}$ will be called the 
\emph{removed ends} of $\widetilde{N^\ast}$, and they are of 
two kinds: those that project onto the internal cusps of $N$ 
are genuine horoballs 
completely contained in $\widetilde{N}$, while a removed end
that project onto a boundary cusp of $N$ is properly contained
in a horoball centered at a point that belongs to the closure at infinity
of two 
different components of $\partial\widetilde{N}$.

For later purposes we insist that $N^\ast\subset N$ is chosen in such
a way that the following conditions hold: the distance between two
distinct removed ends of $\widetilde{N^\ast}$ is at least $1$;
a constant $c>0$ exists such that
the distance between any removed end of $\widetilde{N^\ast}$ 
which projects onto an internal cusp of $N$ and the boundary of 
$\widetilde{N}$
is equal to $c$; a constant $a>0$ exists such that the boundary of
any removed end of $\widetilde{N^\ast}$ which projects onto a 
boundary cusp of $N$ is isometric to $\matR^{n-2}\times [0,a]$.
We point out that the last two conditions imply that every isometry
of $\widetilde{N}$ restricts to an isometry of $\widetilde{N^\ast}$.

We observe that $\widetilde{N^\ast}$ is no longer convex. Moreover,
its boundary is partitioned into the following sets:
the \emph{geodesic boundary} 
$\partial_g \widetilde{N^\ast}=
\partial \widetilde{N^\ast}\cap \partial\widetilde{N}$, which is given
by the union of  
portions of geodesic hyperplanes, and the  \emph{horospherical boundary}
$\partial_h \widetilde{N^\ast}=
\partial \widetilde{N^\ast}\setminus \partial \widetilde{N}$,
which is given by the union of portions of horospheres.

Recall that $d^{\matH}$ is the hyperbolic distance of $\matH^n$
(which restricts of course to a distance on $\widetilde{N^\ast}$,
still denoted by $d^\matH$)
and let
$d^\ast$ be the path distance induced on $\widetilde{N^\ast}$. 
Of course we have $d^{\matH}\leqslant d^\ast$, but
since $\widetilde{N^\ast}$
is not convex, $d^{\matH}$ and $d^\ast$ are in fact quite different
from each other.
In general, $(\widetilde{N^\ast},d^\ast)$ is not quasi-isometric
to $(\widetilde{N^\ast},d^{\matH})$. However, it is easily seen
that $d^{\matH}$ and $d^\ast$ are biLipschitz equivalent below any
given scale,~\emph{i.e.} for any $c>0$ there exists $b\geqslant 1$ 
such that 
\[
b^{-1}\cdot d^{\matH} (x,y)\leqslant d^\ast (x,y)\leqslant
b \cdot d^{\matH} (x,y)
\textrm{ for every}\ x,y\in\widetilde{N^\ast}\ \textrm{with}
\ d^\ast (x,y)\leqslant c.
\]

\paragraph{Quasi-isometries preserve the horospherical boundary}
Let $N_1,N_2$ be hyperbolic $n$-manifolds with non-empty 
geodesic boundary. From now on we also suppose that
if $n=3$, then $\partial N_1$
and $\partial N_2$ are both compact. 
Let $N^\ast_i$ be a compact core of $N_i$ as above, let 
$\widetilde{N}_i\subset\matH^n$ 
(resp.~$\widetilde{N_i^\ast}\subset\matH^n$)
be the universal covering (resp.~the neutered
universal covering) of $N_i$, and $\Gamma_i<{\rm Isom} (\widetilde{N}_i)
<{\rm Isom} (\matH^n)$ be the fundamental group of $N_i$.
We fix a finite set of generators for $\Gamma_i$ and
suppose
that $\varphi:\Gamma_1\to\Gamma_2$ is a $k_1$-quasi-isometry
with $k_1$-pseudo-inverse $\varphi^{-1}$.

Since $N^\ast_i\subset N_i$ is compact,
after fixing  
basepoints $\widetilde{x}_i\in\widetilde{N}_i$, we are provided
with
$k_2$-quasi-isometries $\varphi^\ast:\widetilde{N_1^\ast}\to
\widetilde{N_2^\ast}$, $(\varphi^{-1})^\ast:\widetilde{N_2^\ast}\to
\widetilde{N_1^\ast}$, one the $k_2$-pseudo-inverse of the other,
where $\widetilde{N_i^\ast}$ is endowed
with the path metric $d_i^\ast$, and
$k_2$ only depends on
$k_1$, on the fixed sets of generators for $\Gamma_1$, $\Gamma_2$,
and on the diameters of
$N^\ast_1,N^\ast_2$. 

Recall that our definition of neutered universal covering
implies that a constant $a_i>0$ exists such that
any component of
$\partial_h \widetilde{N_i^\ast}$ is isometric either to
$\matR^{n-1}$ or to $\matR^{n-2}\times [0,a_i]$.
Moreover, if $n=3$ then the boundaries of $N_1$ and of $N_2$
are supposed to be compact. Thus 
if $H$ is a component of
$\partial_h \widetilde{N_i^\ast}$, then $H$ is $r$-quasi-isometric
to $\matR^k$ for some $k\geqslant 2$, where $r$ only depends on $a_i$. 
Therefore an easy application of~\cite[Lemma 3.1]{Sch} 
implies the following:

\begin{lemma}\label{horo:preserved:lemma}
If $H_1$ is a component of 
$\partial_h \widetilde{N_1^\ast}$, then a unique component $H_2$ 
of $\partial_h \widetilde{N_1^\ast}$ exists such that $\varphi^\ast
(H_1)$ is contained in the $r'$-neighbourhood
of $H_2$, where $r'$ is a positive number
depending solely on the geometry of $N^\ast_1,N^\ast_2$ and on $k_2$.
\end{lemma}

Recall that a map that stays at an uniformly bounded distance from
a quasi-isometry is still a quasi-isometry. Moreover, up to swapping the
indices, the statement of Lemma~\ref{horo:preserved:lemma}
also holds for $(\varphi^{-1})^\ast$.
Thus, up to increasing $k_2$ by an amount which depends solely on
$N_1,N_2$ and $k_2$ itself, we can suppose that
$\varphi^\ast$ (resp.~$(\varphi^{-1})^\ast$) 
takes $\partial_h \widetilde{N^\ast_1}$ (resp.~$\partial_h 
\widetilde{N^\ast_2}$)
into $\partial_h \widetilde{N^\ast_2}$ (resp. into~$\partial_h 
\widetilde{N^\ast_1}$).

\paragraph{Extending $\varphi^\ast$ to $\widetilde{N}_1$}

We now describe how
$\varphi^\ast$ can be extended to a quasi-isometry $\widetilde{\varphi}$
from $\widetilde{N}_1$ to
$\widetilde{N}_2$. Even if the construction of $\widetilde{\varphi}$
is essentially
the same as in~\cite[Section 5]{Sch}, we outline it here, since 
some explicit properties 
of $\widetilde{\varphi}$ will be needed in the following paragraphs.
We begin with the following:

\begin{lemma}\label{equivalent:lemma}
Let $(X,d_X), (Y,d_Y)$ be path metric spaces, let $k$ be a positive
constant and suppose the maps
$f:X\to Y$, $g:Y\to X$ have the following properties:
\[
\begin{array}{l}
d_Y (f(x),f(x'))\leqslant k\ \textrm{for every}\ x,x'\in X
\ \textrm{with}\ d_X (x,x')\leqslant 1;\\
d_X (g(y),g(y'))\leqslant k\ \textrm{for every}\ y,y'\in Y
\ \textrm{with}\ d_Y (y,y')\leqslant 1;\\
d_X (g(f(x)),x)\leqslant k\ \textrm{for every}\ x\in X;\\
d_Y (f(g(y)),y)\leqslant k\ \textrm{for every}\ y\in Y.
\end{array}
\]
Then $f$ and $g$ are $\max \{k,3\}$-quasi-isometries.
\end{lemma}
\noindent\emph{Proof:}
Let $x,x'$ be points in $X$ and let $N\in\matN$ be such
that $N\leqslant d_X (x,x')
<N+1$. Since $X$ is a path metric space,
there exist points $x_0=x,x_1,\ldots,x_{N+1}=x'$ with
$d_X (x_i,x_{i+1})\leqslant 1$ for every $i=0,\ldots,N$.
Thus 
\[
d_Y (f(x),f(x'))\leqslant\sum_{i=0}^N d_Y (f(x_i),f(x_{i+1}))
\leqslant Nk+k\leqslant k d_X (x,x')+k.
\] 
The same argument applied to $g$ shows that  
$d_X (g(y),g(y'))\leqslant k d_Y (y,y')+k$ for every $y,y'\in Y$.

Moreover, for every $x,x'\in X$ we have 
\begin{equation*}
\begin{array}{lll}
d_X (x,x') & \leqslant & d_X (x,g(f(x)))+d_X (g(f(x)),g(f(x')))+d_X 
(g(f(x')),x')\\
& \leqslant & k+(k d_Y (f(x),f(x'))+k)+k,
\end{array}
\end{equation*}
whence $k^{-1} d_X (x,x')-3\leqslant d_Y (f(x),f(x'))$.
This shows that $f$ is a $\max \{k,3\}$-quasi-isometric embedding, and the 
conclusion follows
since we can apply the same argument to $g$. 
\finedimo

Before going on, we fix some notation. For any $q\in\partial\matH^n,
x\in\matH^n$ we denote by
$H^q_x$ the horosphere of $\matH^n$ centered at $q$ and containing $x$,
and by $d^q_x$ the natural Euclidean path metric defined on $H^q_x$.
Moreover, we denote by $\xi^q_x$ the projection of $\matH^n$ onto
$H^q_x$, \emph{i.e.}~the map which takes any point $y\in\matH^n$
to the intersection of $H^q_x$ with the geodesic
line containing $y$ and having one endpoint at $q$.
It is easily seen that a universal constant $c_1>0$ exists such that
\[
d_x^q (x,\xi^q_x (x'))\leqslant c_1\ \textrm{for every}\ x,x'\in 
\matH^n\ \textrm{with}\ d^\matH (x,x')\leqslant 1.
\]

Suppose $x$ is a point in $\widetilde{N}_i\setminus\widetilde{N_i^\ast}$.
Then $x$ belongs to exactly one removed end $O\subset \matH^n$. If 
$q\in\partial \matH^n$ is the center of $O$, we denote by $\eta_i (x)$
the intersection of $\partial O$ with the geodesic line
in $\matH^n$ containing $x$ and having an endpoint at $q$
(so $\eta_i (x)=\xi^q_y (x)$, where $y$ is any point belonging to
$\partial O$).

We now define $\widetilde{\varphi}:\widetilde{N}_1\to\widetilde{N}_2$
as follows: if $x$ belongs to $\widetilde{N_i^\ast}$, then
$\widetilde{\varphi} (x)=\varphi^\ast (x)$; otherwise, we set
$\widetilde{\varphi}(x)=y$, where $y$ is the unique point of 
$\widetilde{N}_2\setminus\widetilde{N_2^\ast}$ with 
$\eta_2 (y)=\varphi^\ast (\eta_1 (x))$ and $d^\matH (y,\eta_2 (y))=
d^\matH  (x,\eta_1 (x))$. We say that $\widetilde{\varphi}$ is the
\emph{conical extension} of $\varphi^\ast$, and we denote by
$\widetilde{\varphi}^{-1}$ the conical extension of
$(\varphi^{-1})^\ast$.

\begin{prop}\label{extension:prop}
The map $\widetilde{\varphi}:\widetilde{N}_1\to\widetilde{N}_2$
just constructed is a $k_3$-quasi-isometry, where $k_3$
only depends on $k_2$.
\end{prop}
\noindent\emph{Proof:}
Recall first that 
$d_i^\ast$ and the restriction of $d^\matH$ to $\widetilde{N_i^\ast}$
are biLipschitz equivalent below any given scale for $i=1,2$, 
so a constant $c_2>0$ exists such that
\begin{equation}\label{zero:eq}
d_i^\ast (y,y')\leqslant c_2\ \textrm{for every}\ y,y'\in 
\widetilde{N^\ast_i}\ \textrm{with}\ d^\matH (y,y')\leqslant 2.
\end{equation}

Let $x,x'$ be points in $\widetilde{N}_1$ with 
$d^\matH (x,x')\leqslant 1$. Since 
any two distinct removed ends of $\widetilde{N}_1$ lie at distance
at least $1$ from each other, $x$ and $x'$ do not belong to different
components of $\widetilde{N}_1\setminus \widetilde{N^\ast_1}$.

Suppose first that $x,x'\in\widetilde{N^\ast_1}$. Then 
inequality~(\ref{zero:eq})
implies
\begin{equation}\label{prima:eq}
d^\matH (\widetilde{\varphi} (x),\widetilde{\varphi} (x'))
\leqslant d_2^\ast (\widetilde{\varphi} (x),\widetilde{\varphi} (x'))=
d_2^\ast ({\varphi}^\ast (x),{\varphi}^\ast (x'))\leqslant
k_2 c_2 +k_2.
\end{equation}

\begin{figure}
\begin{center}
\input{ext.pstex_t}
\caption{\small{Extending $\varphi^\ast$ to the removed ends.}}
\label{horo:fig}
\end{center}\end{figure}
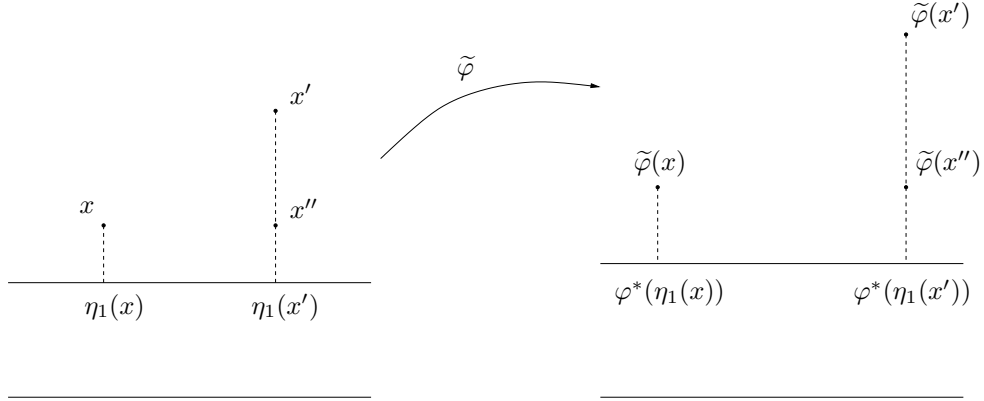

Suppose now that $x$ and $x'$ belong to the same removed end
of $\widetilde{N_1^\ast}$ centered at $q$
and set $b=d^\matH (x,\eta_1 (x))\leqslant
d^\matH (x',\eta_1 (x'))$ (see Figure~\ref{horo:fig}).
Let also $x''=\xi^q_x (x')$ be the projection of $x'$ onto
$H^q_x$. Of course we have $d^\matH (x',\eta_1 (x'))-d^\matH (x,\eta_1 (x))
\leqslant 1$ and $d^q_x (x,x'')\leqslant c_1$.
It is easily seen that $d_1^\ast (\eta_1 (x),\eta_1 (x'))=d^q_{\eta_1 (x)}
(\eta_1 (x),\eta_1 (x'))\leqslant c_1 \exp (b)$, whence
$d_2^\ast (\varphi^\ast (\eta_1 (x)),\varphi^\ast (\eta_1 (x')))\leqslant
k_2 c_1 \exp (b)+k_2$.
Observe now that 
$\widetilde{\varphi} (x),\widetilde{\varphi} (x'),\widetilde{\varphi} (x'')$ belong
to the same removed end of $\widetilde{N^\ast_2}$. Suppose
this end is centered at $s\in\partial\matH^n$. 
Then by construction we have
\[
\begin{array}{lllllll}
d^\matH (\widetilde{\varphi} (x),\widetilde{\varphi} (x')) & \leqslant &
d^s_{\widetilde{\varphi} (x)} (\widetilde{\varphi} (x),
\widetilde{\varphi} (x''))
&=&\exp (-b) d_2^\ast 
(\varphi^\ast (\eta_1 (x)),\varphi^\ast (\eta_1 (x')))\\
&\leqslant&
k_2 c_1 +k_2\exp (-b) &\leqslant & k_2 (c_1+1),
\end{array}
\]
and
\[
d^\matH (\widetilde{\varphi} (x''),\widetilde{\varphi} (x'))
= d^\matH (x'',x') \leqslant 1,
\]
whence
\[
\begin{array}{lll}
d^\matH (\widetilde{\varphi} (x),\widetilde{\varphi} (x'))  \leqslant 
d^\matH (\widetilde{\varphi} (x),\widetilde{\varphi} (x''))+
d^\matH (\widetilde{\varphi} (x''),\widetilde{\varphi} (x'))
\leqslant
k_2 (c_1+1)+1.
\end{array}
\]

Finally suppose that $x\notin\widetilde{N^\ast_1}, x'\in\widetilde{N^\ast_1}$
and set $x''=\eta_1 (x)$.
Of course we have 
$d^\matH (x,x'')\leqslant 1$,
whence $d^\matH (x'',x')\leqslant 2$
by the triangular inequality. 
By construction we have
$d^\matH (\widetilde{\varphi} (x),\widetilde{\varphi} (x''))=
d^\matH (x,x'')\leqslant
1$. Moreover, the same argument yielding 
inequality~(\ref{prima:eq}) also shows that
$d^\matH (\widetilde{\varphi} (x''),\widetilde{\varphi} (x'))\leqslant
k_2 c_2+k_2$. This readily implies 

%\begin{equation}\label{terza:eq}
\[
d^\matH (\widetilde{\varphi} (x),\widetilde{\varphi} (x'))\leqslant
d^\matH (\widetilde{\varphi} (x),\widetilde{\varphi} (x''))+
d^\matH (\widetilde{\varphi} (x''),\widetilde{\varphi} (x))\leqslant
k_2 (c_2+1)+1.
\]
%\end{equation}

Of course, analogous inequalities also hold for $\widetilde{\varphi}^{-1}$.
Moreover, by the very construction of $\widetilde{\varphi},
\widetilde{\varphi}^{-1}$
it follows that for all $x\in\widetilde{N}_1,
y\in\widetilde{N}_2$ we have
\[
\begin{array}{l}
d^\matH (x,\widetilde{\varphi}^{-1}(\widetilde{\varphi} (x)))
\leqslant d^\matH (\eta_1 (x),\widetilde{\varphi}^{-1}
(\widetilde{\varphi} (\eta_1 (x))))\leqslant 
d_1^\ast (\eta_1 (x),\widetilde{\varphi}^{-1}
(\widetilde{\varphi} (\eta_1 (x))))\leqslant k_2,\\
d^\matH (y,\widetilde{\varphi}(\widetilde{\varphi}^{-1} (y)))
\leqslant d^\matH (\eta_2 (y),\widetilde{\varphi}
(\widetilde{\varphi}^{-1} (\eta_2 (y))))\leqslant 
d_2^\ast (\eta_2 (y),\widetilde{\varphi}
(\widetilde{\varphi}^{-1} (\eta_2 (y))))\leqslant k_2.
\end{array}
\]
Now the conclusion follows from Lemma~\ref{equivalent:lemma}.
\finedimo

\paragraph{Quasi-isometries preserve the geodesic boundary}
If $(X,d)$ is a metric space, a 
$k$-quasi-isometric embedding $\alpha:\matR\to X$ is called
\emph{$k$-quasi-geodesic} in $X$. 
The following proposition establishes 
that a quasi-geodesic in hyperbolic space
stays within a finite distance from a genuine geodesic
(see \emph{e.g.}~\cite{BenPet,Ghy} for a proof). 

\begin{prop}\label{quasigeod:prop}
Let $\alpha:\matR\to \matH^n$ be a $k$-quasi-geodesic.
Then there exist a positive number $r$, depending solely
on $k$, and a geodesic $\beta:\matR\to\matH^n$ such that
the image of $\alpha$ is contained in the $r$-neighbourhood
of the image of $\beta$. 
\end{prop}

A quasi-isometry between convex subsets
of hyperbolic space naturally defines a continuous map
between their closures at infinity~\cite{Efr,Tukia2}:

\begin{prop}\label{tukia:prop}
For $i=1,2$ let
Let $W^i$ be a convex subset of $\matH^n$,
let $W^i_\infty\subset \partial\matH^n$ 
be the boundary at infinity of $W^i$ and suppose that
$f:W^1\to W^2$ is a quasi-isometry.
Then there exists a unique homeomorphism 
$\partial f:W^1_\infty\to W^2_\infty$
having the following property:
if $\alpha$ is (the image of) a quasi-geodesic
in $W^1$ with endpoints $x,x'\in W^1_\infty$,
then $f (\alpha)$ is (the image of) a quasi-geodesic
in $W^2_\infty$ with endpoints $\partial f (x),
\partial f (x')\in W^2_\infty$. Such a homeomorphism is called the
extension of $f$ to the conformal boundary of $W^1$.
\end{prop}

We are now ready to prove the following:    

\begin{lemma}\label{boundary:preserved:lemma}
If $S_1$ is a component of 
$\partial \widetilde{N}_1$, then a component $S_2$ 
of $\widetilde{N}_2$ exists such that $\widetilde{\varphi}
(S_1)$ is contained in the $r$-neighbourhood of $S_2$,
where $r$ is a positive number
depending solely on $k_3$.
\end{lemma}
\noindent\emph{Proof:}
Let $\partial \widetilde{\varphi}:\Lambda (\Gamma_1)\to
\Lambda (\Gamma_2)$ be the extension of 
$\widetilde{\varphi}$ to the conformal boundary of $\widetilde{N}_1$.
Since $\partial\widetilde{\varphi}$ 
is a homeomorphism, Lemma~\ref{connect:lemma}
implies that $\partial\widetilde{\varphi}$ takes the closure at
infinity of $S_1$ onto the closure at infinity of a
connected component $S_2$ of $\partial \widetilde{N}_2$.
Let $\alpha$ be (the image of) a geodesic lying on $S_1$.
By Proposition~\ref{quasigeod:prop}, the set
$\widetilde{\varphi} (\alpha)$ is contained in the $r$-neighbourhood
of a geodesic $\beta$ of $\matH^n$. Moreover, $\partial\widetilde{\varphi}$
takes the endpoints of $\alpha$ onto the endpoints of $\beta$,
so $\beta$ lies on $S_2$, and $\widetilde{\varphi} (\alpha)$ is contained 
in the $r$-neighbourhood of $S_2$.
\finedimo

Thus, up to increasing $k_3$ by an amount which depends solely on
$N_1,N_2$ and $k_3$ itself, we can modify $\widetilde{\varphi},
\widetilde{\varphi}^{-1}$ in such a way that
the following conditions hold:
\begin{itemize}
\item
$\widetilde{\varphi}$, $\widetilde{\varphi}^{-1}$
are $k_3$-quasi-isometries between 
$(\widetilde{N}_1,d^\matH)$ and $(\widetilde{N}_2,d^\matH)$, which are
one the $k_3$-pseudo-inverse of the other;
\item
$\widetilde{\varphi}$ (resp.~$\widetilde{\varphi}^{-1}$)
takes 
$\widetilde{N_1^\ast}$ (resp.~$\widetilde{N_2^\ast}$) into 
$\widetilde{N_2^\ast}$ (resp.~$\widetilde{N_1^\ast}$);
\item
$\widetilde{\varphi}$ (resp.~$\widetilde{\varphi}^{-1}$)
takes 
$\partial_h \widetilde{N_1^\ast}$ (resp.~$\partial_h \widetilde{N_2^\ast}$) into 
$\partial_h \widetilde{N_2^\ast}$ (resp.~$\partial_h \widetilde{N_1^\ast}$);
\item
$\widetilde{\varphi}$, $\widetilde{\varphi}^{-1}$
restrict to $k_3$-quasi-isometries between 
$(\widetilde{N^\ast_1},d_1^\ast)$ and $(\widetilde{N^\ast_2},d_2^\ast)$, which are
one the $k_3$-pseudo-inverse of the other;
\item
$\widetilde{\varphi}$ (resp.~$\widetilde{\varphi}^{-1}$)
is the conical extension of
$\widetilde{\varphi}|_{N_1^\ast}$ (resp.~of $\widetilde{\varphi}^{-1}
|_{N_2^\ast}$);
\item
$\widetilde{\varphi}$ (resp.~$\widetilde{\varphi}^{-1}$)
takes the geodesic boundary of
$\widetilde{N}_1$ (resp.~of $\widetilde{N}_2$) into the geodesic boundary of
$\widetilde{N}_2$ (resp.~of $\widetilde{N}_1$).
\end{itemize}

Also observe that if $S,S'$ are distinct hyperplanes of $\matH^n$,
then $S$ is not contained in any $r$-neighbourhood of $S'$, so  
 $\widetilde{\varphi}$
induces a bijection between the components of
$\partial \widetilde{N}_1$ and the components of
$\partial \widetilde{N}_2$.

\paragraph{Mirroring along the boundary}  
Let now $DN_i$ be the manifold obtained by mirroring $N_i$ along its
boundary. Since $\partial N_i$ is totally geodesic, $DN_i$ is a hyperbolic
$n$-manifold without boundary. We can choose a universal 
covering $p_i :\matH^n\to DN_i$ in such a way that
$p_i$ restricts to the given universal covering $\widetilde{N}_i\to
N_i$. The preimages under $p_i$ of $N_i$ and of its mirror copy 
define a tessellation $\calT_i$ of $\matH^n$ whose pieces are convex polyhedra
isometric to $\widetilde{N}_i$. Since $\widetilde{\varphi}$ induces
a bijection between the components of 
$\partial \widetilde{N}_1$ and the components of
$\partial \widetilde{N}_2$, the tessellations $\calT_1$ 
and $\calT_2$ are combinatorially equivalent to each other.
Thus the map $\widetilde{\varphi}$ can be extended to the whole
of $\matH^n$ just by repeatedly mirroring along the hyperplanes
which bound the pieces of $\calT_1$ and of $\calT_2$. More precisely,
 unique maps $\overline{\varphi},\overline{\varphi}^{-1}:\matH^n\to\matH^n$
exist with the following properties:
\begin{itemize}
\item
$\overline{\varphi} (\widetilde{N}_1)=\widetilde{N}_2$,
$\overline{\varphi}^{-1} (\widetilde{N}_2)=\widetilde{N}_1$;
\item
$\overline{\varphi} |_{\widetilde{N}_1}=\widetilde{\varphi}$,
$\overline{\varphi}^{-1} |_{\widetilde{N}_2}=\widetilde{\varphi}^{-1}$;
\item
the image under $\overline{\varphi}$ (resp. under~$\overline{\varphi}^{-1}$)
of a piece of $\calT_1$ (resp.~of $\calT_2$) is contained in
a piece of 
$\calT_2$ (resp.~of $\calT_1$);
\item
adjacent pieces of $\calT_1$ (resp.~of $\calT_2$) are taken by
$\overline{\varphi}$ (resp.~by $\overline{\varphi}^{-1}$) into adjacent pieces of 
$\calT_2$ (resp.~of $\calT_1$);
\item
let $P_1,P_1'$ be adjacent pieces of $\calT_1$ with $P_1\cap P_1'=S_1$,
and set $\overline{\varphi}(P_1)=P_2$, $\overline{\varphi} (P'_1)=P'_2$, $S_2=P_2\cap P'_2$.
If $\sigma_i$ is the hyperbolic reflection along $S_i$, then
we have $(\overline{\varphi}\circ\sigma_1)|_{P_1}=(\sigma_2\circ\overline{\varphi})|_{P_1}$,
$(\overline{\varphi}^{-1}\circ\sigma_2)|_{P_2}=(\sigma_1\circ\overline{\varphi}^{-1}|_{P_2}$.
\end{itemize}

Let $\Theta_i$ be the preimage under $p_i$ of the double of  
$N^\ast_i$. Then $\Theta_i$ is obtained from
$\matH^n$ by removing a countable set of horoballs at distance
at least one from each other. 
Observe that for $i,j=1,2$, if $f:\Theta_i \to \Theta_j$ is a
map taking the boundary of $\Theta_i$ into the boundary of
$\Theta_j$, then it does make sense to speak of the 
conical extension of $f$, which is defined on the whole of $\matH^n$.
We denote by $\calT^\ast_i$ the tessellation
of $\Theta_i$ whose pieces are obtained by intersecting the pieces
of $\calT_i$ with $\Theta_i$. Every piece of $\calT^\ast_i$
is the image of $\widetilde{N^\ast_i}$ under an element of
${\rm Isom}(\matH^n)$.
Let $\delta^\ast_i$ be the path metric on $\Theta_i$ induced by the 
hyperbolic Riemannian structure. If $P^\ast_i$ is a piece of $\calT_i^\ast$,
then $(P^\ast_i,\delta^\ast_i |_{P^\ast_i\times P^\ast_i})$ is isometric to
$(\widetilde{N^\ast_i},d_i^\ast)$. Moreover, the restriction
of $\overline{\varphi}$ to any piece $P^\ast_1$ of $\calT_1^\ast$ defines a $k_3$-quasi-isometry
between $P^\ast_1$ and the piece of $\calT^\ast_2$ containing 
$\overline{\varphi} (P^\ast_1)$, both endowed with their path distance.

\begin{lemma}\label{double:lemma}
The maps $\overline{\varphi},\overline{\varphi}^{-1}$ are $k_4$-quasi-isometries
of $\matH^n$, where $k_4$ depends solely on $k_3$.
\end{lemma}
\noindent\emph{Proof:}
We first prove that $\overline{\varphi},\overline{\varphi}^{-1}$ restrict to quasi-isometries between
$(\Theta_1,\delta_1^\ast)$ and $(\Theta_2,\delta_2^\ast)$ 
that are one the pseudo-inverse of the other.
Since $N^\ast_i$ is compact, 
it is easily seen
that a positive constant $c$ exists such that if $x,x'$ belong to 
different connected components of the geodesic boundary
of a piece $P^\ast_i$ of $\calT_i^\ast$, then $\delta_i^\ast (x,x')\geqslant
d^\matH (x,x')\geqslant c$. Let us consider points $x,x'\in\Theta_1$
with $cN\leqslant \delta_1^\ast (x,x')<c(N+1)$. Since
$\Theta_1$ is a geodesic metric space, 
there exist an integer $N'\leqslant N$ and 
points $x_0=x,x_1,\ldots,x_{N'\! +2}=x'$ with 
the following properties:  
$\sum_{i=0}^{N'\! +1} \delta_1^\ast (x_i,x_{i+1})=\delta_1^\ast (x,x')$,
and $x_i,x_{i+1}$ belong to the same closed piece of
$\calT_1^\ast$ for every $i=0,\ldots,N'\! +1$.
Since $\overline{\varphi}$ restricts to a $k_3$-quasi-isometry on each piece
of $\calT_1^\ast$, we obtain 
\[
\begin{array}{lllll}
\delta_2^\ast (\overline{\varphi}(x),\overline{\varphi}(x'))&\leqslant&
\sum_{i=0}^{N'\! +1} 
\delta_2^\ast (\overline{\varphi}(x_i),\overline{\varphi}(x_{i+1}))
&\leqslant &\sum_{i=0}^{N'\! +1} (k_3 \delta_1^\ast (x_i,x_{i+1})+k_3)\\
& \leqslant & k_3 \delta_1^\ast (x,x') +k_3 (N'\! +2) & \leqslant
& k_3 (1+1/c)\delta_1^\ast (x,x') +2k_3.
\end{array}
\]
The same argument also shows that $\delta_1^\ast 
(\overline{\varphi}^{-1}(y),\overline{\varphi}^{-1}
(y'))\leqslant  k_3 (1+1/c)\delta_2^\ast (y,y')+ 2k_3$ for all $y,y'\in
\Theta_2$. Suppose now that $x$ belongs to a piece $P^\ast_1$ of $\calT_1^\ast$
and let $P^\ast_2$ be the piece of $\calT_2^\ast$ containing $\overline{\varphi} (P^\ast_1)$.
It is easily seen that $\overline{\varphi}^{-1}(P^\ast_2)$ is contained in $P^\ast_1$.
Moreover, by the very construction of $\overline{\varphi}$ 
it follows that elements $g_1,g_2$
of ${\rm Isom}(\matH^n)$ exist such that $g_i (\widetilde{N_i})=P^\ast_i$,
$\overline{\varphi}|_{P^\ast_1}=g_2 \widetilde{\varphi} g_1^{-1}|_{P^\ast_1}$
and $\overline{\varphi}^{-1}|_{P^\ast_2}=g_1\widetilde{\varphi}^{-1} g_2^{-1}|_{P^\ast_2}$.
This easily implies that 
$\delta_1^\ast (x,\overline{\varphi}^{-1}(\overline{\varphi}(x)))=
\delta_1^\ast (g_1^{-1} (x), \widetilde{\varphi}^{-1}
(\widetilde{\varphi}(g_1^{-1}(x))))\leqslant k_3$. The same argument also
shows that $\delta_2^\ast 
(y,\overline{\varphi}(\overline{\varphi}^{-1}(y)))\leqslant k_3$
for every $y\in\Theta_2$. By Lemma~\ref{equivalent:lemma}, this implies
that $\overline{\varphi},\overline{\varphi}^{-1}$ restrict to quasi-isometries between
$(\Theta_1,\delta_1^\ast)$ and $(\Theta_2,\delta_2^\ast)$ 
that are one the pseudo-inverse of the other.

Since $\overline{\varphi},\overline{\varphi}^{-1}$ are the conical extensions
of their restrictions to $\Theta_1,\Theta_2$, the conclusion
follows from
the same argument 
yielding Proposition~\ref{extension:prop}. 
\finedimo

\paragraph{Quasi-conformal homeomorphisms}
Let an isometric identification of $\matH^n$ with the Poincar\'e disc
model of hyperbolic $n$-space 
be fixed from now on, and put on $\partial\matH^n=S^{n-1}$
the Riemannian structure induced by the usual Euclidean metric on
$\matR^n$. We observe that such a structure is compatible with the 
canonical conformal structure of $\partial\matH^n$, and we denote 
by $d_\partial$ the induced path metric on $\partial\matH^n$.
For $q\in\partial\matH^n,\, \epsilon>0$ we set
$B(q,\epsilon)=\{x\in\partial\matH^n:\, d_\partial (q,x)\leqslant \epsilon\}$.
Let $f:\partial\matH^n\to \partial\matH^n$ be a homeomorphism.
For every $q\in \partial\matH^n$ we denote by $K(f)_q\in [1,\infty]$ the number
\[
K(f)_q=\limsup_{\epsilon\to 0} \frac{\sup \{d_\partial (f(x),f(q)):
\, x\in\partial B(q,\epsilon)\} }{\inf \{d_\partial (f(x),f(q)):
\, x\in\partial B(q,\epsilon)\} }.
\]
The value $K(f)_q$ gives a measure of how far $f$ is from being conformal
at $q$,
and does not depend on the choice of the 
identification of $\matH^n$ with the Poincar\'e model of hyperbolic $n$-space.
We say that $f$ is $K$-quasi-conformal at $q$ if $K(f)_q\leqslant K$, and that
$f$ is $K$-quasi-conformal if it is $K$-quasi-conformal at every point
of $\partial\matH^n$. Moreover, $f$ is said to be quasi-conformal
if it is $K$-quasi-conformal for some $K<\infty$.

Suppose $f$ is differentiable at $q$ and $df_q$ is invertible, 
and recall that
the Riemannian structure on $\partial\matH^n$ defines a metric on 
$T(\partial\matH^n)_q$ and on $T(\partial\matH^n)_{f(q)}$. If 
$C\subset T(\partial\matH^n)_q$ is a round sphere (not necessarily
centered in $0$), then $df_q(C)\subset
T(\partial\matH^n)_{f(q)}$ is an ellipsoid. We denote by 
$K' (f)_q$ the ratio between the longest and the shortest axes of
$df_q(C)$, and observe that $K'(f)_q$ is indeed well-defined, \emph{i.e.}~it
does not depend on the choice of the sphere $C$. 
The following lemma is straightforward:

\begin{lemma}\label{linear:lemma}
Let $f:\partial \matH^n\to\partial \matH^n$ be a 
homeomorphism, let 
$q\in\partial\matH^n$ and suppose that
the differential of 
$f$ at $q$ exists and is invertible. Then
$K' (f)_q=K (f)_q$.
\end{lemma}

The following fundamental results are taken from~\cite{Mos}:

\begin{prop}\label{mostowbis:prop}
Let $f:\partial \matH^n\to\partial \matH^n$ be a 
quasi-conformal homeomorphism. Then the differential of 
$f$ exists and is invertible  almost
everywhere (with respect to the Lebesgue measure). 
\end{prop}

\begin{prop}\label{mostow:prop}
Suppose $f:\partial \matH^n\to\partial \matH^n$ is a 
quasi-conformal homeomorphism that is $1$-quasi-conformal almost
everywhere (with respect to the Lebesgue measure). Then
$f$ is differentiable and $1$-quasi-conformal everywhere, so it is
the trace at infinity of a unique isometry of $\matH^n$.
\end{prop}

\paragraph{Extending $\overline{\varphi}$ to the conformal boundary}
By Proposition~\ref{tukia:prop}, the quasi-isometries
$\overline{\varphi},\overline{\varphi}^{-1}$ uniquely define extensions
$\partial{\overline{\varphi}},\partial{\overline{\varphi}}^{-1}:
\partial{\matH}^n\to\partial{\matH}^n$ that are one the inverse of
the other (see 
also~\cite{Mos2,thu,BenPet}).  
Moreover, $\partial \overline{\varphi},\partial \overline{\varphi}^{-1}$ are
quasi-conformal homeomorphisms~\cite{Mos2}.

Recall now that $p_i:\matH^n\to DN_i$ is the universal covering
that extends the universal covering $\widetilde{N}_i\to N_i$. 
Let $D\Gamma_i<{\rm Isom}(\matH^n)$ by the group
of automorphisms of the covering $p_i$. We have of course
$DN_i\cong\matH^n/ D\Gamma_i$ and $\Gamma_i<D\Gamma_i$.
A geodesic hyperplane in $\matH^n$ is a \emph{face}
of $\calT_i$ if it is a boundary component of some piece of
$\calT_i$.
We set
\[
\begin{array}{cll}
\calC_i &=&\{C\subset\partial\matH^n:\ C\ \textrm{is the boundary at infinity
of a face of}\ \calT_i\},\\
D\Lambda_i &=& D\Gamma_i \cdot \Lambda (\Gamma_i)\subset\partial\matH^n,
\end{array}
\] 
and we observe that $\bigcup_{C\in\calC_i} C\subset D\Lambda_i$. Moreover,
if $C\in\calC_1$ then $\partial\overline{\varphi} (C)\in\calC_2$.
It is well-known that the limit set of any geometrically finite subgroup
of ${\rm Isom}(\matH^n)$ has either null or full Lebesgue measure~\cite{Ahl}
(recently, the same statement has been proved to hold for any discrete 
finitely generated subgroup of ${\rm Isom}(\matH^n)$~\cite{Agol,Cal}).
Thus $D\Lambda_i$ is a countable union of measure zero subsets of
$\partial\matH^n$, and has itself zero Lebesgue measure.

\paragraph{Deforming $\overline{\varphi}$ into an isometry}
Recall that we have fixed on $\partial\matH^n$ a Riemannian metric 
induced by an identification of $\matH^n$ with the Poincar\'e disc
model of hyperbolic $n$-space. 
We say that a point $q\in\partial\matH^n$ is 
\emph{secluded} (with respect to $\calT_1$)
if a sequence $\{B^i\}_{i\in\matN}$ of round closed 
balls in $\partial \matH^n$
exists with the following properties: $q$ belongs to the interior
of $B^i$ for every 
$i\in\matN$, the diameter of $B^i$ tends to $0$ as $i$ tends to infinity,
and the boundary of $B^i$ is an element of $\calC_1$
for every $i\in\matN$. 

\begin{lemma}\label{secluded:lemma}
Let $q$ be a point in $\partial\matH^n\setminus D\Lambda_1$. Then
$q$ is secluded with respect to $\calT_1$.
\end{lemma}
\noindent\emph{Proof:}
Let $\alpha:[0,\infty)\to\matH^n$ be any geodesic ray with endpoint $q$.
If $\alpha$ intersects a finite number of faces of $\calT_1$, then
$q$ is contained in the closure at infinity of a piece of $\calT_1$,
whence in $D\Lambda_1$, a contradiction.
Thus there exists an infinite set $\{F^i\}_{i\in\matN}$ of faces
of $\calT_1$ such that
$F^i$ intersects the 
geodesic $\alpha$ in a single point $q^i$
and $\lim_{i\to\infty} q^i
=q$ in $\overline{\matH}^n$.
For $i\in\matN$ let $C^i$ be the boundary at infinity of
$F^i$. 
Since $q\notin D\Lambda_1$, for every $i\in\matN$ 
we have $q\notin C^i$, so a well-defined round ball $B^i\subset
\partial \matH^n$ exists whose boundary is equal to $C^i$ and whose
interior contains $q$. Also observe that (up to passing to a subsequence)
we can assume that $C^{i+1}$ is contained in $B^i$ for every 
$i\in\matN$.
In order to prove the lemma we only have 
to show that the diameter of $B^i$ tends to $0$ as $i$ tends to infinity.

Suppose this is not true. In this case,
since $C^{i+1}\subset B^{i}$, it is easily seen that a limit round
circle $C^{\infty}\subset \partial\matH^n$ exists such that
$\lim_{i\to\infty} C^i=C^{\infty}$ in the Haussdorff topology
on closed subsets of $\partial\matH^n$.
Let $F^{\infty}\subset \partial\matH^n$ be the hyperbolic hyperplane
bounded by $C^{\infty}$. Then $\lim_{i\to\infty} 
F^i=F^{\infty}$ in the Haussdorff topology
on closed subsets of $\overline{\matH}^n$. Since 
the union of the faces of $\calT_1$ is closed in
${\matH}^n$, this easily implies that $F^{\infty}$ is a face
of $\calT_1$, so $C^{\infty}$ belongs to $\calC_1$ and
is contained in $D\Lambda_1$. Since $\lim_{i\to\infty} q^i=q$
we should have $q\in\overline{F}^\infty\cap
\partial\matH^n = C^{\infty}$, whence $q\in D\Lambda_1$, a contradiction.
\finedimo

\begin{lemma}\label{conformal:lemma}
Let $q\in\partial\matH^n$ be secluded with respect to $\calT_1$
and suppose that $d(\partial\overline{\varphi})_q$ exists and is invertible. 
Then $\partial\overline{\varphi}$ is $1$-quasi-conformal at $q$.
\end{lemma}
\noindent\emph{Proof:}
Let us identify $\partial \matH^n$ with $\matR^{n-1}\cup\{\infty\}$
in such a way that $q$ corresponds to $0\in\matR^{n-1}$.
Without loss of generality, we can also assume that $\partial \overline{\varphi} (q)=0
\in\matR^{n-1}\subset\partial\matH^n$.
Then there exists a family
$\{B^i\}_{i\in\matN}$ of round 
balls in $\matR^{n-1}\subset \partial \matH^n$ with the following properties:
$q\in{\rm int}\, B^i$ for every 
$i\in\matN$; the Euclidean diameter of $B^i$
tends to $0$ as $i$ tends to infinity; 
the boundary of $B^i$ is an element $C^i$ of $\calC_1$
for every $i\in\matN$. 

Let
$\lambda_i:\matR^{n-1}\to\matR^{n-1}$ be the linear dilation having ratio
equal to the Euclidean diameter of $B^i$. If we identify $T (\partial\matH^n)_q$ and
$T (\partial \matH^n)_{\partial\overline{\varphi} (q)}$ with
$T (\matR^{n-1})_0=\matR^{n-1}$, then 
\begin{equation}\label{lim:eq}
\lim_{i\to\infty} \lambda_i^{-1} \circ
\partial \overline{\varphi} \circ {\lambda_i}=
d (\partial\overline{\varphi})_q\quad\textrm{uniformely on any compact subset of}
\ \matR^{n-1}.
\end{equation}
Let $B(0,1)\subset \matR^{n-1}$ be the closed unitary ball centered in $0$,
and observe that
the round sphere $L^i={(\lambda_i)^{-1}}(C^i)$ has Euclidean diameter
equal to $1$ and is contained in 
$B(0,1)$. Since $\partial\overline{\varphi}$ sends spheres in $\calC_1$ onto
spheres in $\calC_2$,  the sphere $L^i$ is
sent by $\lambda_i^{-1} \circ
\partial \overline{\varphi} \circ \lambda_i$ onto a round sphere. 
On the other hand, the liner map $d (\partial\overline{\varphi})_q$ sends
$L^i$ onto an ellipsoid $E^i$. 
Let $\nu_i$ be the ratio between the 
longest and the shortest axes of $E^i$.  
Since
$d (\partial\overline{\varphi})_q$ is invertible,
there exists $\delta>0$ such that the shortest axis of $E^i$ 
is longer than $\delta$ for every $i\in\matN$.
Together with equality~(\ref{lim:eq}), this implies that
$\lim_{i\to\infty} \nu_i=1$. But $\nu_i=K' (f)_q$ for every $i$,
whence $K' (f)_q=1$. Now the conclusion follows from  
Lemma~\ref{linear:lemma}.
\finedimo

We recall that $D\Lambda_i$ has zero Lebesgue measure, so by 
Proposition~\ref{mostowbis:prop}, Lemma~\ref{secluded:lemma}, 
Lemma~\ref{conformal:lemma} and Proposition~\ref{mostow:prop}
the map $\partial\overline{\varphi}$ is the trace at infinity
of an element $\widehat{\varphi}$ of ${\rm Isom} (\matH^n)$.

\begin{lemma}\label{boundeddist:lemma}
We have $d^\matH (\overline{\varphi} (x),\widehat{\varphi} (x))\leqslant r$ for every
$x\in\matH^n$, where $r$ is
a positive constant depending solely on $k_4$.
\end{lemma}
\noindent\emph{Proof:}
Let $\alpha_1,\alpha_2\subset\matH^n$ be two geodesic lines
intersecting perpendicularly at $x$. 
By Proposition~\ref{quasigeod:prop}
the point $\overline{\varphi} (x)$
belongs both to the $r'$-neighbourhood of 
$\widehat{\varphi}(\alpha_1)$
and to the $r'$-neighbourhood of
$\widehat{\varphi} (\alpha_2)$, where $r'$ only depends on $k_4$. 
But $\widehat{\varphi}(\alpha_1)$
intersects 
$\widehat{\varphi} (\alpha_2)$ perpendicularly at
$\widehat{\varphi} (x)$, so $d^\matH (\overline{\varphi} (x),
\widehat{\varphi} (x))\leqslant r$, where $r$ only depends on $r'$, whence
on $k_4$.
\finedimo

By construction, the trace at infinity of $\widehat{\varphi}$ sends
the boundary at infinity of any component of $\partial \widetilde{N}_1$
into the boundary at infinity of a component of $\partial \widetilde{N}_2$.
Since $\widetilde{N}_i$ is the hyperbolic convex hull of its geodesic
boundary for $i=1,2$, this implies that 
$\widehat{\varphi}$ restricts to an isometry
between $\widetilde{N}_1$ and $\widetilde{N}_2$.

\section{Main results}\label{corollari:section}
For $i=1,2$ let $N_i$ be a hyperbolic $n$-manifold with non-empty geodesic boundary,
and suppose that if $n=3$, then $\partial N_i$ is compact. Let $\widetilde{N}_i\subset
\matH^n$ be the universal covering of $N_i$ and $\Gamma_i<{\rm Isom} (\widetilde{N}_i)
<{\rm Isom} (\matH^n)$ be the fundamental group of $N_i$. 
Fix a finite set $S_i$ of generators of $\Gamma_i$, fix 
a basepoint $\widetilde{x}_i\in\widetilde{N}_i$, and consider $\Gamma_i$ 
as a subset
of $\widetilde{N}_i$ via the embedding $\Gamma_i\ni\gamma\mapsto\gamma 
(\widetilde{x}_i)\in\widetilde{N}_i$.
The following statement summarizes the results obtained in the preceding section.

\begin{teo}\label{fundamental:teo}
For any $k\geqslant 1$ there exists $r=r(k)>0$ such that if 
$\varphi: \calC (\Gamma_1,S_1)\to 
\calC (\Gamma_2,S_2)$ is a $k$-quasi-isometry, then an
isometry $\widehat{\varphi}:\widetilde{N}_1\to\widetilde{N}_2$ exists
such that $d^\matH (\widehat{\varphi} (\gamma), \varphi (\gamma))\leqslant r$
for every $\gamma\in\Gamma_1\subset\widetilde{N}_1$.
\end{teo}

\paragraph{Quasi-isometry implies commensurability}
We can now easily prove Theorem~\ref{main1:teo}.
Let $N_1, N_2$ be as above and suppose that $\pi_1 (N_1)$
is quasi-isometric to $\pi_1 (N_2)$. By Theorem~\ref{fundamental:teo},
$N_1$ and $N_2$ have isometric universal coverings.
By Lemma~\ref{important:lemma}, this implies that $N_1$ is commensurable with $N_2$.

\paragraph{Commensurator and quasi-isometry group}
We now prove Theorem~\ref{comm1:teo}.
Let $N$ be a hyperbolic $n$-manifold with non-empty 
geodesic boundary, let $\widetilde{N}\subset\matH^n$ 
(resp.~$\widetilde{N^\ast}\subset\matH^n$) 
be the universal covering (resp.~the neutered universal covering)
of $N$, and denote by $d^\ast$ the path metric on $\widetilde{N^\ast}$.
As before, we fix a point $\widetilde{x}\in\widetilde{N^\ast}\subset\widetilde{N}$, and
consider the fundamental group 
$\Gamma<{\rm Isom} (\widetilde{N})<{\rm Isom} (\matH^n)$ of $N$ as a subset
of $\widetilde{N^\ast}\subset\widetilde{N}$.

Let $\varphi$ be a quasi-isometry of $\Gamma$ into itself, and denote by $[\varphi]$
the equivalence class of $\varphi$ in ${\rm QIsom} (\Gamma)$. By 
Theorem~\ref{fundamental:teo}, there exist $r>0$ and an isometry
$\widehat{\varphi}:\widetilde{N}\to\widetilde{N}$ such that
$d^\matH (\varphi (\gamma),\widehat{\varphi}(\gamma))\leqslant r$ for every
$\gamma\in\Gamma$. It is easily seen that this condition uniquely determines
$\widehat{\varphi}$, and that $\widehat{\varphi}$ 
only depends on the equivalence class
of $\varphi$. Thus the map 
\[
{\rho}:{\rm QIsom} (\Gamma)\to
{\rm Isom}(\widetilde{N}),\quad \rho ([\varphi])=\widehat{\varphi}
\]
gives a well-defined group homomorphism.

Injectivity of $\rho$ is obvious. Moreover, 
any isometry of $\widetilde{N}$ restricts to an isometry
of $(\widetilde{N^\ast},d^\ast)$, whence to a quasi-isometry of $\Gamma\subset
\widetilde{N^\ast}$. This implies that $\rho$ is surjective.
Recall now that the left translation by an element $\gamma\in\Gamma$ 
naturally defines
an element of ${\rm QIsom} (\Gamma)$, still denoted by $\gamma$. 
By construction we have 
$d^\matH
(\gamma(\widetilde{x}),
\rho(\gamma)(\widetilde{x}))\leqslant r$ for every
$\gamma\in\Gamma$, and this easily implies that $\rho(\gamma)=\gamma$
for every
$\gamma\in\Gamma$. This concludes the proof
of Theorem~\ref{comm1:teo}.

\paragraph{Quasi-isometric rigidity}
We now prove Theorem~\ref{main3:teo}. We keep the notation
from the preceding
paragraph. Let $G$ be a finitely generated group and suppose
$\varphi:G\to \Gamma=\pi_1 (N)$ is a $k$-quasi-isometry
with $k$-pseudo-inverse $\varphi^{-1}$.
Without loss of generality, we can also suppose that
$\varphi (1_G)=1_\Gamma$ and $\varphi^{-1} (1_\Gamma) = 1_G$. 
Fix an element $g$ in $G$. The left translation by $g$ defines
an isometry of any Cayley graph of $G$, whence a $k'$-quasi-isometry
$\mu (g)$ of $\Gamma$, where $k'$ depends on $k$ but not on $g$.
By Theorem~\ref{fundamental:teo},
there exist a constant $r>0$ and an isometry
$\nu (g)\in{\rm Isom}(\widetilde{N})$ such that
$d^\matH (\mu (g)(\gamma)(\widetilde{x}), \nu (g) 
(\gamma(\widetilde{x})))\leqslant r$
for every $\gamma\in \Gamma$. Taking $\gamma=1_\Gamma$,
we get in particular 
\begin{equation}\label{discrete:eq}  
d^\matH (\varphi (g) (\widetilde{x}),\nu (g)(\widetilde{x}))\leqslant r,
\end{equation}
where $r$ does not depend on $g$.

It is easily seen that the map $\nu:G\to{\rm Isom} (\widetilde{N})$
is a group homomorphism. We are going to show that 
$\nu$ has finite kernel and discrete image. Let $\{g_n\}_{n\in\matN}\subset G$
be any sequence such that $\lim_{n\to\infty} \nu (g_n)=1_\Gamma$.
Up to discarding a finite number of terms, by equation~(\ref{discrete:eq})
we can suppose that 
\begin{equation}\label{discrete2:eq}  
d^\matH (\varphi (g_n)(\widetilde{x}),  
\widetilde{x})\leqslant r+1\quad\textrm{for every}\ n\in\matN.
\end{equation}
Since $\Gamma$ acts properly discontinuously on $\widetilde{N}$, 
this implies that the set 
$\{\varphi (g_n)\}_{n\in\matN}
\subset \Gamma$ is finite, and since
$\varphi:G\to\Gamma$ is a quasi-isometry, this gives in turn that
$\{g_n\}_{n\in\matN}\subset G$ is a finite set. Thus $\nu (G)
\subset{\rm Isom}(\widetilde{N})$ is discrete and $\nu$ has finite kernel.

Let now $q$ be a point in 
$\Lambda (\Gamma)=\Lambda ({\rm Isom}(\widetilde{N}))$
and let $\{\gamma_n\}_{n\in\matN}\subset\Gamma$ be a sequence with
$\lim_{n\to\infty} \gamma_n (\widetilde{x})=q$.  
Since $\varphi$ is a $k$-quasi-isometry and $d^\matH\leqslant
d^\ast$ on $\widetilde{N^\ast}$, 
we have $d^\matH (\varphi(\varphi^{-1}(\gamma_n))
(\widetilde{x}), \gamma_n (\widetilde{x}))\leqslant r'$, where
$r'$ is a positive constant which does not depend on $n$. 
Thus, if $g_n=\varphi^{-1} (\gamma_n)$, by equation~(\ref{discrete:eq})
we have 
\[
d^\matH (\nu (g_n) (\widetilde{x}),\gamma_n (\widetilde{x})) \leqslant 
d^\matH (\nu (g_n) 
(\widetilde{x}),\varphi(g_n)(\widetilde{x}))
+ d^\matH (\varphi(g_n)(\widetilde{x}), \gamma_n (\widetilde{x}))
\leqslant  r'',
\]
where $r''$ is a positive constant which does not depend on $n$.
Thus $\lim_{n\to\infty} \nu (g_n) (\widetilde{x})=
\lim_{n\to\infty} \gamma_n (\widetilde{x})=q$, and
$\Lambda (\nu (G))=\Lambda (\Gamma)$. Recall now that if $\Delta$ is a geometrically
finite subgroup of ${\rm Isom} (\matH^n)$ whose limit set is not equal to
the whole of $\partial\matH^n$, then any finitely generated subgroup
of $\Delta$ is geometrically finite~(see 
\emph{e.g.}~\cite[Proposition 7.1]{Mor}). In our context, this implies
that $\nu (G)$ is geometrically finite, so the convex core
$N_G:=\widetilde{N}/\nu (G)$ is a finite-volume orbifold.
If $N_{\rm orb}=\widetilde{N}/{\rm Isom} (\widetilde{N})$, then
\[
[{\rm Isom} (\widetilde{N}):\nu (G)]=
{\rm vol}\, N_G/{\rm vol}\, N_{\rm orb}<\infty.
\]
This shows that $\nu (G)$ is commensurable with $\Gamma$.

\section{Counterexamples}\label{contro:section} 
  
We now show that the conclusions of Theorems~\ref{main1:teo},
\ref{comm1:teo}
are no longer true if we consider hyperbolic
$3$-manifolds with non-compact geodesic boundary.

\paragraph{Non-commensurable manifolds sharing the same 
fundamental group}

We begin by 
recalling that Thurston's hyperbolization theorem
for Haken manifolds~\cite{thu2} gives necessary and sufficient topological
conditions on a $3$-manifold to be hyperbolic with
geodesic boundary:

\begin{teo}\label{hyper:teo}
Let $\overline{M}$ be a compact orientable $3$-manifold 
with non-empty boundary, let $\calT$ be the set of
boundary tori of $\overline{M}$ and
let $\calA$ be a family of disjoint closed annuli in 
$\partial \overline{M}\setminus\calT$. 
Then $M=\overline{M}\setminus (\calT\cup\calA)$ is hyperbolic
if and only if the following conditions hold:
\begin{itemize}
\item
the components of $\partial M$ 
have negative Euler characteristic;
\item
$\overline{M}\setminus\calA$ 
is boundary-irreducible and geometrically atoroidal;
\item
the only proper essential annuli contained in $M$ are parallel
in $\overline{M}$ to the annuli in $\calA$.
\end{itemize}
\end{teo}

The following proposition readily implies Theorem~\ref{contro1:teo}.

\begin{prop}\label{contro:prop}
Let $N$ be any orientable hyperbolic $3$-manifold with compact non-empty
geodesic boundary.
Then there exists a hyperbolic $3$-manifold with non-compact 
geodesic boundary
which is homotopically equivalent to $N$ but not commensurable with $N$.
\end{prop}
\noindent\emph{Proof:}
Let $\alpha$ be a simple essential loop in
$\partial N$. We define $N'$ as $N\setminus\alpha$ and note that $N$ 
and $N'$ have a common compactification $\overline{N}=\overline{N'}$
such that $N=\overline{N}\setminus\calT$,
$N'=\overline{N}\setminus(\calT\cup\calA)$,
where $\calT$ is the family of the boundary tori of $\overline{N}$ 
and $\calA$ is a closed regular neighbourhood of $\alpha$ in
$\partial N$.
Moreover, it is easily seen that since $(\overline{N},
\calT,\emptyset)$ satisfies the assumptions of Theorem~\ref{hyper:teo},
so does $(\overline{N},
\calT,\calA)$, so $N'$ is hyperbolic. 
Of course $N'$ is homotopically equivalent to $N$.

Since $\partial N$ is compact, there do not exist different
components of $\partial \widetilde{N}\subset \matH^3$ 
that meet in $\partial \matH^n$.
On the other hand, since $\partial N'$ is non-compact a pair
$(S^1,S^2)$ of components of $\partial \widetilde{N'}\subset \matH^3$  
exists such that $\overline{S}^{1}\cap\overline{S}^{2}=\{ {\rm pt.}\}\in
\partial\matH^3$. Thus $\widetilde{N}$ is not isometric to $\widetilde{N'}$,
and $N$ is not commensurable with $N'$.
\finedimo  

\paragraph{A hyperbolic $3$-manifold with free non-Abelian
fundamental group}
In this paragraph we prove Theorem~\ref{contro2:teo}.
To this aim we construct a hyperbolic $3$-manifold $M$
with non-compact geodesic boundary
which is homotopically
equivalent to a genus-$2$ handlebody. The fundamental group of 
$M$ is isomorphic to $\matZ\ast\matZ$, and it is well-known
that (equivalence classes of) left translations define a subgroup
of \emph{infinite} index in the quasi-isometry group of
$\matZ\ast\matZ$. 
On the other hand, by Lemma~\ref{fare:lemma} the group
${\rm Comm} (\pi_1 (M))/\pi_1 (M)$ is finite, so $M$ provides
the manifold required in the statement of  Theorem~\ref{contro2:teo}.

The following
explicit construction of $M$ is taken from~\cite{Fri}.
\begin{figure}
\begin{center}
\input{prima6.pstex_t}
\caption{\small{The manifold $M$  
is obtained by gluing in pairs the non-shadowed
faces of the regular ideal octahedron 
along suitable isometries.}}\label{octa:fig}
\end{center}\end{figure}
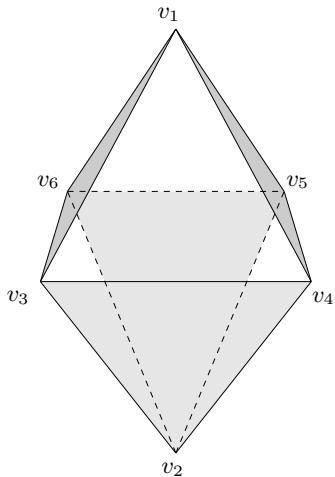
Let $O\subset\matH^3$ be the regular ideal octahedron and let
$v_1,\ldots,v_6$ be the vertices of $O$ as shown in Figure~\ref{octa:fig}.
We denote by $F_{ijk}$ the face of $O$ with vertices 
$v_i,v_j,v_k$. Let $g:F_{134}\to F_{156}$ be the unique 
isometry such that $g(v_1)=v_1$, $g(v_3)=v_6$ and $g(v_4)=v_5$, and 
$h:F_{236}\to F_{254}$ be the unique 
isometry such that $h(v_2)=v_2$, $h(v_3)=v_4$,
$h(v_6)=v_5$.
We define $M$ to be the manifold obtained by gluing $O$
along $g$ and $h$.   
Since all the dihedral angles of $O$ are right, it is 
easily seen that the metric on $O$ induces a complete
finite-volume hyperbolic structure on $M$ 
such that the shadowed 
faces in Figure~\ref{octa:fig} are glued along their edges to give
a non-compact totally geodesic boundary.

\begin{figure}
\begin{center}
\input{handle.pstex_t}
\caption{\small{The natural compactification of $M$ is a
genus-$2$ handlebody with a boundary annulus. 
Here the annulus is represented by 
its core curve.}}\label{handle:fig}
\end{center}
\end{figure}
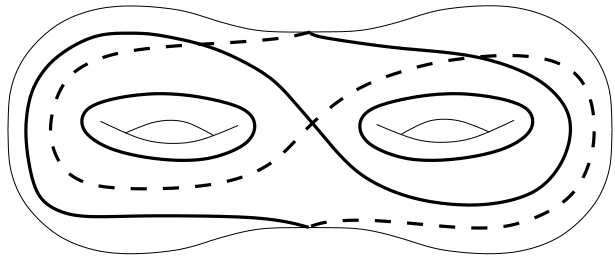

From a topological and combinatorial point of view, 
an ideal octahedron with four marked faces as in Figure~\ref{octa:fig}
is equivalent to a truncated tetrahedron with the edges connecting
truncation triangles removed, which is
\begin{figure}
\begin{center}
\input{tetra.pstex_t}
\caption{\small{An ideal octahedron, a truncated tetrahedron, 
and a tetrapod
with arcs.}}\label{tetra:fig}
\end{center}
\end{figure}
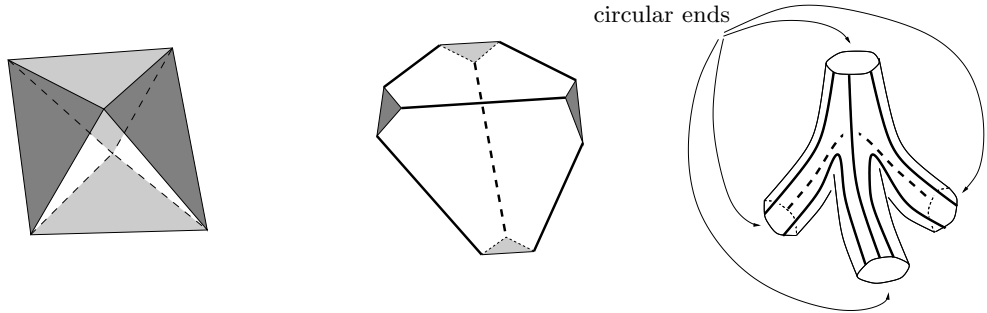
in turn equivalent to a ``tetrapod'' with six arcs connecting
circular ends removed, as shown in
Figure~\ref{tetra:fig}. Under this identification, the four shadowed
ideal faces of $O$ correspond to the four regions into which the 
lateral surface of the tetrapod is cut by the $6$ arcs, while
the non-shadowed ideal faces of $O$
correspond to the four discs at the ends of the tetrapod. Therefore 
the manifold $M$  
is obtained from the tetrapod by suitably gluing together in pairs
the discs at its four ends. So this manifold
is homeomorphic to a handlebody with boundary 
loops removed.
Using this correspondence we can easily draw a picture of
the natural compactification of $M$.
This picture is shown in Figure~\ref{handle:fig}.
For a more detailed description of the natural compactification 
of hyperbolic $3$-manifolds with non-compact geodesic boundary
obtained by 
gluing regular ideal octahedra see~\cite{CFMP}.

\bibliographystyle{amsalpha}
\bibliography{bibliocomm}

\vspace{1.5 cm}

\noindent
\hspace*{6cm}Dipartimento di Matematica L.~Tonelli\\ 
\hspace*{6cm}largo Bruno Pontecorvo, 5 \\
\hspace*{6cm}56127 Pisa, Italy\\ 
\hspace*{6cm}frigerio@mail.dm.unipi.it
\vspace{.5 cm}

%\noindent
%\hspace*{6cm}Dipartimento di Matematica\\ 
%\hspace*{6cm}Universit\`a di Pisa\\ 
%\hspace*{6cm}Via F. Buonarroti 2\\ 
%\hspace*{6cm}56127 Pisa, Italy\\ 
%\hspace*{6cm}martelli@mail.dm.unipi.it
%\vspace{.5 cm} 

%\noindent
%\hspace*{6cm}Dipartimento di Matematica Applicata\\ 
%\hspace*{6cm}Universit\`a di Pisa\\
%\hspace*{6cm}Via Bonanno Pisano 25B\\ 
%\hspace*{6cm}56126 Pisa, Italy\\ 
%\hspace*{6cm}petronio@dm.unipi.it

\end{document}

%% file: defs.tex
\newtheorem{lemma}{Lemma}[section] 
\newtheorem{teo}[lemma]{Theorem}
 
\newtheorem{prop}[lemma]{Proposition}
\newtheorem{cor}[lemma]{Corollary}

\newcommand{\matN}{\ensuremath {\mathbb{N}}}
\newcommand{\matR} {\ensuremath {\mathbb{R}}}

\newcommand{\matZ} {\ensuremath {\mathbb{Z}}}

\newcommand{\matH} {\ensuremath {\mathbb{H}}}

\newcommand{\calC} {\ensuremath {\mathcal{C}}}
\newcommand{\calA} {\ensuremath {\mathcal{A}}}

\newcommand{\calT} {\ensuremath {\mathcal{T}}}

\newcommand{\nota} [1] {\caption{\footnotesize{#1}}}

\newfont{\Got}{eufm10 scaled 1200}

\font\titsc=cmcsc10 scaled 1200

\newcommand{\finedimo}{{\hfill\hbox{$\square$}\vspace{2pt}}}

%% file: ext.pstex_t
\begin{picture}(0,0)%
\special{psfile=ext.pstex}%
\end{picture}%
\setlength{\unitlength}{1579sp}%
\begingroup\makeatletter\ifx\SetFigFont\undefined
% extract first six characters in \fmtname
\def\x#1#2#3#4#5#6#7\relax{\def\x{#1#2#3#4#5#6}}%
\expandafter\x\fmtname xxxxxx\relax \def\y{splain}%
\ifx\x\y   % LaTeX or SliTeX?
\gdef\SetFigFont#1#2#3{%
  \ifnum #1<17\tiny\else \ifnum #1<20\small\else
  \ifnum #1<24\normalsize\else \ifnum #1<29\large\else
  \ifnum #1<34\Large\else \ifnum #1<41\LARGE\else
     \huge\fi\fi\fi\fi\fi\fi
  \csname #3\endcsname}%
\else
\gdef\SetFigFont#1#2#3{\begingroup
  \count@#1\relax \ifnum 25<\count@\count@25\fi
  \def\x{\endgroup\@setsize\SetFigFont{#2pt}}%
  \expandafter\x
    \csname \romannumeral\the\count@ pt\expandafter\endcsname
    \csname @\romannumeral\the\count@ pt\endcsname
  \csname #3\endcsname}%
\fi
\fi\endgroup
\begin{picture}(15024,6225)(589,-6373)
\put(1801,-5011){\makebox(0,0)[lb]{\smash{\SetFigFont{10}{12.0}{rm}$\eta_1 (x)$}}}
\put(4426,-5011){\makebox(0,0)[lb]{\smash{\SetFigFont{10}{12.0}{rm}$\eta_1 (x')$}}}
\put(5026,-1711){\makebox(0,0)[lb]{\smash{\SetFigFont{10}{12.0}{rm}$x'$}}}
\put(7651,-1261){\makebox(0,0)[lb]{\smash{\SetFigFont{10}{12.0}{rm}$\widetilde{\varphi}$}}}
\put(1726,-3436){\makebox(0,0)[lb]{\smash{\SetFigFont{10}{12.0}{rm}$x$}}}
\put(5026,-3511){\makebox(0,0)[lb]{\smash{\SetFigFont{10}{12.0}{rm}$x''$}}}
\put(10426,-2761){\makebox(0,0)[lb]{\smash{\SetFigFont{10}{12.0}{rm}$\widetilde{\varphi} (x)$}}}
\put(14851,-2761){\makebox(0,0)[lb]{\smash{\SetFigFont{10}{12.0}{rm}$\widetilde{\varphi} (x'')$}}}
\put(14776,-436){\makebox(0,0)[lb]{\smash{\SetFigFont{10}{12.0}{rm}$\widetilde{\varphi} (x')$}}}
\put(13876,-4786){\makebox(0,0)[lb]{\smash{\SetFigFont{10}{12.0}{rm}$\varphi^\ast (\eta_1 (x'))$}}}
\put(10126,-4786){\makebox(0,0)[lb]{\smash{\SetFigFont{10}{12.0}{rm}$\varphi^\ast(\eta_1 (x))$}}}
\end{picture}

%% file: prima6.pstex_t
\begin{picture}(0,0)%
\special{psfile=prima6.pstex}%
\end{picture}%
\setlength{\unitlength}{2368sp}%
\begingroup\makeatletter\ifx\SetFigFont\undefined
% extract first six characters in \fmtname
\def\x#1#2#3#4#5#6#7\relax{\def\x{#1#2#3#4#5#6}}%
\expandafter\x\fmtname xxxxxx\relax \def\y{splain}%
\ifx\x\y   % LaTeX or SliTeX?
\gdef\SetFigFont#1#2#3{%
  \ifnum #1<17\tiny\else \ifnum #1<20\small\else
  \ifnum #1<24\normalsize\else \ifnum #1<29\large\else
  \ifnum #1<34\Large\else \ifnum #1<41\LARGE\else
     \huge\fi\fi\fi\fi\fi\fi
  \csname #3\endcsname}%
\else
\gdef\SetFigFont#1#2#3{\begingroup
  \count@#1\relax \ifnum 25<\count@\count@25\fi
  \def\x{\endgroup\@setsize\SetFigFont{#2pt}}%
  \expandafter\x
    \csname \romannumeral\the\count@ pt\expandafter\endcsname
    \csname @\romannumeral\the\count@ pt\endcsname
  \csname #3\endcsname}%
\fi
\fi\endgroup
\begin{picture}(3195,4980)(3196,-5461)
\put(6121,-2446){\makebox(0,0)[lb]{\smash{\SetFigFont{9}{10.8}{\updefault}$v_5$}}}
\put(6391,-3661){\makebox(0,0)[lb]{\smash{\SetFigFont{9}{10.8}{\updefault}$v_4$}}}
\put(4816,-5461){\makebox(0,0)[lb]{\smash{\SetFigFont{9}{10.8}{\updefault}$v_2$}}}
\put(4771,-691){\makebox(0,0)[lb]{\smash{\SetFigFont{9}{10.8}{\updefault}$v_1$}}}
\put(3511,-2446){\makebox(0,0)[lb]{\smash{\SetFigFont{9}{10.8}{\updefault}$v_6$}}}
\put(3196,-3661){\makebox(0,0)[lb]{\smash{\SetFigFont{9}{10.8}{\updefault}$v_3$}}}
\end{picture}

%% file: tetra.pstex_t
\begin{picture}(0,0)%
\special{psfile=tetra.pstex}%
\end{picture}%
\setlength{\unitlength}{1026sp}%
\begingroup\makeatletter\ifx\SetFigFont\undefined
% extract first six characters in \fmtname
\def\x#1#2#3#4#5#6#7\relax{\def\x{#1#2#3#4#5#6}}%
\expandafter\x\fmtname xxxxxx\relax \def\y{splain}%
\ifx\x\y   % LaTeX or SliTeX?
\gdef\SetFigFont#1#2#3{%
  \ifnum #1<17\tiny\else \ifnum #1<20\small\else
  \ifnum #1<24\normalsize\else \ifnum #1<29\large\else
  \ifnum #1<34\Large\else \ifnum #1<41\LARGE\else
     \huge\fi\fi\fi\fi\fi\fi
  \csname #3\endcsname}%
\else
\gdef\SetFigFont#1#2#3{\begingroup
  \count@#1\relax \ifnum 25<\count@\count@25\fi
  \def\x{\endgroup\@setsize\SetFigFont{#2pt}}%
  \expandafter\x
    \csname \romannumeral\the\count@ pt\expandafter\endcsname
    \csname @\romannumeral\the\count@ pt\endcsname
  \csname #3\endcsname}%
\fi
\fi\endgroup
\begin{picture}(23595,7558)(6678,-9119)
\put(20851,-2086){\makebox(0,0)[lb]{\smash{\SetFigFont{9}{10.8}{\updefault}${\rm circular\ ends}$}}}
\end{picture}

%% file: handle.pstex_t
\begin{picture}(0,0)%
\special{psfile=handle.pstex}%
\end{picture}%
\setlength{\unitlength}{1776sp}%
\begingroup\makeatletter\ifx\SetFigFont\undefined
% extract first six characters in \fmtname
\def\x#1#2#3#4#5#6#7\relax{\def\x{#1#2#3#4#5#6}}%
\expandafter\x\fmtname xxxxxx\relax \def\y{splain}%
\ifx\x\y   % LaTeX or SliTeX?
\gdef\SetFigFont#1#2#3{%
  \ifnum #1<17\tiny\else \ifnum #1<20\small\else
  \ifnum #1<24\normalsize\else \ifnum #1<29\large\else
  \ifnum #1<34\Large\else \ifnum #1<41\LARGE\else
     \huge\fi\fi\fi\fi\fi\fi
  \csname #3\endcsname}%
\else
\gdef\SetFigFont#1#2#3{\begingroup
  \count@#1\relax \ifnum 25<\count@\count@25\fi
  \def\x{\endgroup\@setsize\SetFigFont{#2pt}}%
  \expandafter\x
    \csname \romannumeral\the\count@ pt\expandafter\endcsname
    \csname @\romannumeral\the\count@ pt\endcsname
  \csname #3\endcsname}%
\fi
\fi\endgroup
\begin{picture}(8454,3487)(589,-4503)
\end{picture}